\input amstex.tex
\documentstyle{amsppt}
\baselineskip 16pt plus 2pt
\advance\vsize -2\baselineskip
\parskip 2pt
\NoBlackBoxes
\input epsf

\input xy
\xyoption{all}
\def\dofig#1#2{\epsfxsize=#1 \centerline{\epsfbox{#2}}}

\topmatter

\title A short course on Witten Helffer-Sj\"ostrand  theory \endtitle

\thanks   Supported in part by NSF
\endthanks

\author D. Burghelea (Ohio State University)
\endauthor
\leftheadtext {Witten Helffer Sj\"ostrand  theory}
\rightheadtext {D. Burghelea}

\abstract

Witten-Helffer-Sj\"ostrand theory is an addition to
Morse theory and Hodge-de Rham theory for Riemannian
manifolds and considerably improves on them by injecting some spectral
theory of elliptic operators.
It can  serve as a general tool to prove results about comparison
of
numerical invariants associated to compact manifolds analytically,
i.e. by using a Riemannian metric, or combinatorially, i.e. by using a
triangulation.
It can be also refined to provide an alternative presentation of
Novikov Morse
theory and improve on it in many respects. In particular it can be
used
in symplectic topology and in dynamics.
This material represents my Notes for a three lectures course given
at the
Goettingen summer school on groups and geometry, June 2000.

\endabstract
\toc
\widestnumber\subhead{3.1}
\head 0. Introduction  \endhead

\head 1. Lecture 1: Morse theory revisited \endhead

       a. Generalized triangulations

       b. Morse Bott generalized triangulations

       c. G-generalized triangulations

       d. Proof of Theorem 1.1

       e. Appendix : The transversality of the maps $p$ and $s$ from section d. 

\head 2. Lecture 2: Witten deformation and the spectral properties of
the Witten Laplacian
\endhead

       a. De Rham theory and integration

       b. Witten deformation and Witten Laplacians

       c. Spectral gap theorems

       d. Applications

       e. Sketch of the proof of Theorem 2.1

\head 3. Lecture 3: Helffer Sj\"ostrand Theorem, on an asymptotic
improvement
of the Hodge-de Rham  theorem\endhead

       a. Hodge de Rham Theorem 

       b. Helffer Sj\"ostrand Theorem (Theorem 3.1)

       c. Extensions and a survey of other applications

\head  References. \endhead

\endtoc

\endtopmatter

\document

\vskip .5in

\proclaim{0. Introduction}\endproclaim

Witten Helffer Sj\"ostrand theory, or abbreviated, WHS -theory,
consists
of a number of results which considerably improve on Morse theory and
De Rham -Hodge theory.

The intuition behind  the WHS -theory is provided  by
physics and consists in regarding a compact smooth  manifold equipped
with a
Riemannian metric
and a Morse function (or closed 1-form) as an interacting system of
harmonic
oscillators. This
intuition was first exploited by E. Witten, cf[Wi], in order to
provide
a short "physicist's proof " of Morse inequalities, a rather simple
but very
useful result in topology.

Helffer and
Sj\"ostrand have completed Witten's
picture with their results on Schr\"odinger operators and
have considerably strengthened
Witten's mathematical statements, cf~\cite {HS2}. The work of Helffer
and
Sj\"ostrand on the Witten theory
can be substantially simplified  by using simple observations more or
less
familiar to topologists, cf [BZ] and [BFKM].

The mathematics behind the WHS-theory is almost entirely based on the
following two
facts: the
existence of a gap in the spectrum of the Witten Laplacians (a one parameter family of 
deformed  Laplace-Beltrami operator involving a Morse function $h$), detected
by elementary mini-max characterization of the spectrum of selfadjoint
positive
operators and simple estimates involving the equations
of the harmonic oscillator.  The
Witten Laplacians in the neighborhood of critical points in ``
admissible
coordinates'' are given by such equations.
The simplification we referred  to are due to a compactification
theorem
for the space of trajectories and of the unstable sets of the gradient
of a
Morse function with respect to a "good" Riemannian metric. This
theorem can be
regarded as a strengthening of the basic results of elementary Morse
theory.

The theory, initially considered for a Morse function, can be easily
extended to a Morse
closed one form, and even to the more general case of a Morse-Bott
form.These are closed 1-forms
which in some neighborhood of a connected component of the zero set
are differential
of a Morse-Bott function.
So far the theory has been very useful to obtain an alternative
derivation of results concerning 
the comparison  of numerical
invariants
associated to compact manifolds analytically (i.e. by using a
Riemannian
metric,) and combinatorially cf[BZ1], [BZ2], [BFKM], [BFK1],
[BFK4], [BFKM], [BH].

The theory provides an alternative (analytic) approach to the Novikov-
Morse theory
with considerable improvements and consequently has applications in
symplectic topology and dynamics. These aspects will be developed in a
forthcoming paper [BH].

This minicourse is presented as a series of three lectures.
The first is a reconsideration of elementary Morse theory, with the
sketch of the
proof of the compactification theorem.
The second discusses the "Witten deformation" of the Laplace Beltrami operator 
and its implications and
provides a sketch of
the proof of Theorem 2.1 the main result of the section.
The third presents Helffer-Sj\"ostrand Theorem as an asymptotic
improvement of the 
Hodge-de Rham Theorem and finally surveys some of the existing applications.
Part of the material presented in these notes will contained in a book
[BFK5] in preparation, which will be
written in
collaboration with L. Friedlander and T. Kappeler.

\vskip .5in

{\bf Lecture 1: Morse Theory revisited. }

\vskip .2in

\subhead {a. Generalized triangulations}
\endsubhead

\vskip.2in
Let $ M^n $ be a compact closed smooth manifold of dimension $ n $.
A {\it generalized triangulation }is provided by a pair $(h,g),$
$ h: M \to \Bbb R $
a smooth function, $ g $ a Riemannian metric so that 

{\bf C1.}  For any critical point $ x $ of $ h $ there exists a
coordinate chart in
the
neighborhood of $ x $ so that
in these  coordinates $ h $ is quadratic  and $ g $ is Euclidean.

More precisely, for any critical point $ x $ of $ h,$ ($x \in Cr(h)$),
there exists a coordinate chart
$\varphi: (U,x)\to (D_{\varepsilon},0),\  U$ an open neighborhood of
$x$ in $M,\
D_{\varepsilon}$ an open disc of
radius $\varepsilon$ in $ \Bbb R^n,\  \varphi$ a diffeomorphism with
$ \varphi(x)=0, $
so that :
$$(i)\   h\cdot \varphi^{-1}(x_1, x_2, \cdots, x_n)=c - 1/2( x_1^2+
\cdots x_k^2) +
1/2 ( x_{k+1}^2+\cdots x_n^2)$$
$$(ii)\  (\varphi^{-1})^*(g)\  \text{is given by}\  g_{ij}(x_1,x_2,
\cdots,x_n)= \delta
_{ij}$$
Coordinates so that (i) and (ii) hold are
called {\it admissible}.

It follows then that any critical point $x\in Cr(h)$ has a well
defined index,
$i(x)=\text{ index}(x)=k,$
$k$ the number
of the negative squares in the expression (i), which is independent of
the
choice of  a coordinate system (with respect to which $h$ has the form
(i)).

Consider the vector field $-grad_g(h)$ and for any
$y\in M,$
denote  by $ \gamma_y(t), -\infty < t < \infty,$ the unique trajectory
of
$-grad_g(h)$ which satisfies
the condition $\gamma_y(0)=y.$
\newline For $x\in Cr(h)$ denote by $ W_x^- $ resp. $ W_x^+ $ the sets
$$ W_x^{\pm}=\{y\in M| \lim_{t\to \pm\infty}\gamma_y(t)=x\}. $$
In view of (i), (ii) and of the theorem of existence,
unicity and smooth
dependence on the initial condition
for the solutions of ordinary differential equations, $ W_x^- $
resp. $ W_x^+ $ is a smooth submanifold diffeomorphic to $ \Bbb R^k $
resp. to $\Bbb R^{n-k},$ with $k=\text{ index} (x).$ This can be
verified easily based on
the
fact that
$$ \varphi(W_x^-\cap U_x)=\{(x_1,x_2,\cdots,x_n)\in D(\varepsilon) |x_
{k+1}=x_{k+2}=
\cdots= x_n=0\},$$
and
$$ \varphi(W_x^+\cap U_x)=\{(x_1,x_2,\cdots,x_n)\in D(\varepsilon) |x_
{1}=x_{2}=
\cdots= x_k=0\}.$$

Since $ M $ is compact and C1 holds, the set $ Cr(h) $ is finite and
since $M$
is closed (i.e. compact and without boundary),
$M= \bigcup_{x\in Cr(h)}W_x^-.$ As
already observed each $W_x^-$
is a smooth submanifold diffeomorphic to
$ \Bbb R^k, $ $ k= $index$(x),$ i.e. an open cell.

{\bf C2.} The vector field $-grad_g h$ satisfies the Morse-Smale
condition if for any
$x,y \in Cr(h)$, $W_x^- $
and $ W_y^+$
 are transversal.

C2 implies that $\Cal M(x,y):= W_x^- \cap W_y^+$ is a smooth manifold
of dimension
equal to
 $\text{index}(x) - \text{index}(y).$ $\Cal M(x,y)$
is equipped with the action $\mu: \Bbb R\times \Cal M(x,y) \to
 \Cal M(x,y),$
defined by $\mu(t, z) = \gamma _z(t).$

If $\text{index}(x) \leq \text{index}(y),$
 and $x\neq y,$ in view of the transversality requested by the Morse
 Smale
condition,
$\Cal M(x,y) = \emptyset.$

If $x\neq y$ and $\Cal M(x,y) \neq \emptyset,$
the action $\mu$ is free and we denote the quotient
$\Cal M(x,y)/\Bbb R $
by $\Cal T(x,y);$
$ \Cal T(x,y)$ is a smooth manifold of dimension 
$\text{index}(x) \leq \text{index}(y)-1,$
diffeomorphic to the submanifold
$ \Cal M(x,y) \cap h^{-1}(\lambda),$ for any real number $\lambda$ in
the open
interval $(\text{index}\ x,  \text{index}\ y).$ The elements of
$\Cal T(x,y)$
are the trajectories from
``$x$ to $y$'' and such a trajectory will usually  be
denoted by $\gamma.$

If $ x=y, $ then $W_x^- \cap W_x^+ = {x}.$

Further the condition C2 implies that the
partition of $M$ into open cells is actually a smooth cell complex.
To formulate this fact precisely we
recall that an
\newline $n-$dimensional  manifold $X$ with corners is a paracompact
Hausdorff space
equipped with a
maximal smooth atlas with charts
$\varphi : U\to \varphi(U)\subseteq \Bbb R^n_+$
with $ \Bbb R^n_+ = \{(x_1, x_2,\cdots x_n) | x_i\geq 0\}.$  The
collection of points
of $ X $
which correspond (by some and then by any chart) to points in
$\Bbb R^n$
with exactly $k$ coordinates equal to zero is
a well defined subset of $X$ and it will be denoted by $X_k.$ It
has a structure of a smooth
$(n-k)-$dimensional manifold.
$\partial X = X_1 \cup X_2 \cup \cdots X_n$
is a closed subset which is a topological manifold and
$(X,\partial X)$
is a topological manifold with boundary $\partial X.$ A compact
smooth manifold
with corners, $X,$ with interior diffeomorphic to the Euclidean space,
will be called a compact smooth cell.

For any string of critical points $x=y_0, y_1, \cdots,  y_k$ with
$$\text{index}(y_0) > \text{index}(y_1)>,\cdots,> \text{index}(y_k),$$
consider the smooth manifold of dimension  $\text{index}\ y_0 -k,$
$$\Cal T(y_0,y_1)\times\cdots \Cal T(y_{k-1},y_k)\times W_{y_k}^-,$$
and the smooth
map $$i_{y_0, y_1,\cdots ,y_k}: \Cal T(y_0,y_1)\times\cdots
\times \Cal T(y_{k-1},y_k)\times W_{y_k}^-
\to M,$$
defined by
$i_{y_0, y_1,\cdots ,y_k}(\gamma_1,\cdots, \gamma_k, y):= i_{y_k}(y)$,
where
$i_x: W_x^- \to M$ denotes the inclusion of
$W_x^-$ in $M.$

\proclaim {Theorem 1.1}
Let $\tau =(h,g)$  be a generalized triangulation.

1)For any critical point $x\in Cr(h)$
the smooth manifold $W_x^-$ has a canonical compactification
$\hat W_x^-$ to a
compact manifold
with corners and the inclusion $i_x$ has a smooth extension
$\hat i_x: \hat W_x^- \to M$ so that :

\noindent (a) $(\hat W_x^-)_k = \bigsqcup_{(x,y_1,\cdots, y_k)}
\Cal T(x,y_1)\times \cdots\times
\Cal T(y_{k-1}, y_k)\times W_{y_k}^-,$

\noindent (b) the restriction of $\hat i_x$ to
$\Cal T(x,y_1)\times \cdots \times
\Cal T(y_{k-1}, y_k)\times W_{y_k}^-$ is given by
\newline $i_{x, y_1\cdots, y_k}.$

2) For any two critical points $x, y$ with $i(x) > i(y) $ the
smooth
manifold $\Cal T(x,y)$ has a canonical compactification
$\hat {\Cal T}(x,y)$ to a compact
manifold
with corners and

$\hat {\Cal T}(x,y)_k = \bigsqcup_{(x,y_1,\cdots, y_k=y)}
\Cal T(x,y_1)\times \cdots\times
\Cal T(y_{k-1}, y_k).$
\endproclaim

The proof of Theorem 1.1 will be sketched in the last subsection of
this section.

This theorem was probably well known to experts before it was
formulated by Floer in
the framework
of  infinite dimensional Morse theory cf.~\cite {F}. As formulated,
Theorem 1.1 is
stated in [AB]. The proof sketched in [AB] is excessively complicated
and incomplete.
A considerably simpler proof will be sketched in Lecture 1 subsection
d) and is contained
in [BFK5] and [BH].
\vskip .2in
{\bf Observation:}

O1: The name of generalized triangulation for $\tau= (h,g)$ is
justified
by the fact that any simplicial
smooth triangulation  can be obtained as a generalized triangulation,
cf~\cite{Po}.

O2: Given a Morse function $h$
and a Riemannian
metric $g,$ one can perform arbitrary small $C^0-$ perturbations of
$g,$
so that the pair consisting of $h$ and the perturbed metric is a
generalized triangulation, cf [Sm] and [BFK5].

Given a generalized triangulation $\tau= (h,g),$ and for any critical
point $x\in Cr(h)$ an orientation
$\Cal O_x$ of $W_x^-,$ one can associate a cochain complex of vector
spaces
over the field $\Bbb K$ of real or complex numbers, $(C^*(M,\tau),
\partial^*).$ Denote the collection of these orientations by $o.$
The differential $\partial^*$ depends on the chosen
orientations $o:= \{\Cal O_x |x\in Cr(h)\}.$
To describe this complex we introduce the incidence numbers
$$ I_q: Cr(h)_q\times Cr(h)_{q-1} \to \Bbb Z$$ defined as follows:
\newline If $\Cal T (x,y)=\emptyset,$ we put $I_q(x,y)= 0.$
\newline If $\Cal T (x,y) \ne \emptyset,$  for any
$\gamma\in \Cal T (x,y),$ the set
$\gamma\times W_y^-$ is an open subset of the boundary
$\partial\hat{W}_x^-$ and the orientation
$\Cal O_x$ induces an orientation on it.
If this is the same as the orientation $\Cal O_y,$
we set  $\varepsilon(\gamma)= +1,$ otherwise we set
$\varepsilon(\gamma)= -1.$ Define $I_q(x,y)$ by
$$I_q(x,y)= \sum_{\gamma
\in \Cal T (x,y)} \varepsilon(\gamma).$$
In the case $M$ is an oriented manifold, the orientation of
$M$ and the orientation $\Cal O_x$ on $W_x^-$ induce
an orientation $\Cal O_x^+$ on the stable
manifold $W_x^+.$ 

For any $c$ in the open interval $(h(y), h(x)),$
$h^{-1}(c)$ carries a canonical orientation induced from the
orientation of
$M.$
One can check  that $I_q(x,y)$ is the intersection number of
$W_x^-\cap h^{-1}(c)$ with $W_y^+\cap h^{-1}(c)$ inside $h^{-1}(c)$
and is also the incidence number of the open cells
$W_x^-$ and $W_y^-$
in the $CW-$ complex structure provided by $\tau.$

Denote by $(C^*(M,\tau),\partial^*_{(\tau,o)})$ the cochain
complex of $\Bbb K-$ Euclidean vector spaces defined by

\noindent (1)  $C^q(M,\tau):= Maps (Cr_q(h), \Bbb K)$

\noindent (2)
$\partial^{q-1}_{(\tau,o)}: C^{q-1}(M,\tau) \to
 C^q(M,\tau),
\  (\partial ^{q-1}f) (x)=
\sum _{y\in Cr_{q-1}(h)} I_q(x,y)f(y),$

\noindent where $x\in Cr_q(h).$

\noindent (3) Since $C^q(M,\tau)$ is equipped with a canonical base
provided by the
maps  $E_x$ defined by $E_x(y)= \delta_{x,y},$
$x,y\in Cr_q(h),$
it carries a natural scalar product which makes
$E_x,$ $x\in Cr_q(h),$ orthonormal.

\proclaim {Proposition 1.2}
For any $q, \  \partial^q\cdot \partial ^{q-1}=0.$
\endproclaim
A geometric proof of this Proposition follows from Theorem 1.1.
The reader can also derive it by observing that
$(C^*(M,\tau),\partial^*)$ as defined is nothing but
the cochain complex associated to the $CW-$complex structure
provided by $\tau$  via Theorem 1.1.

\vskip .2in

\subhead {\bf b. Morse Bott generalized triangulations}
\endsubhead

The concept of generalized triangulation and Theorem 1.1 above can be
extended  to pairs
$(h,g)$
with $h$ a Morse Bott function; i.e $Cr(h)$ consists of a disjoint
union of
compact connected smooth submanifolds $\Sigma$ and the Hessian of
$h$ at any $x\in \Sigma$
is nondegenerated
in the
normal directions of $\Sigma.$ More precisely a {\it MB-generalized
triangulation} is a pair $\tau= (h,g)$
which
satisfies C'1 and C'2 below:

{\bf C'1:} $Cr(h)$ is a disjoint union of closed connected
submanifolds
 $\Sigma,$  and for any $\Sigma$ there the exist {\it admissible}
 coordinates in some 
neighborhood of $\Sigma.$ An admissible coordinate chart around $\Sigma$ is
provided by:

1)  two
orthogonal vector
bundles $\nu_{\pm}$ over $\Sigma$ equipped with scalar product
preserving connections
(parallel transports)$\nabla_{\pm}$ so that $\nu_+\oplus \nu_-$
is isomorphic to the
normal bundle of $\Sigma,$

2) a closed tubular neighborhood of $\Sigma,$
$\varphi: (U, \Sigma )\to (D_{\varepsilon}(\nu_-\oplus
 \nu_-),\Sigma),$
$U$ closed neighborhood of $\Sigma,$
so that :

$(i): h\cdot \varphi^{-1}(v_1, v_2)=c - 1/2 ||v_1||^2+ 1/2 ||v_2||^2$
where $v_\pm\in E(\nu_\pm);$

$(ii): (\varphi^{-1})^*(g)$ is the metric induced from the restriction
of $g$ on
$\Sigma,$
the scalar products and the connections in $\nu_{\pm}.$

The rank of $\nu_-$ will be called the index of $\Sigma$ and denoted
by $i(\Sigma)=\text{index}(\Sigma).$
As before for any $x \in Cr(h)$ consider $W^{\pm}_x$ and introduce
$W^{\pm}_{\Sigma}= \cup_{x\in \Sigma}W^{\pm}_x$.

{\bf C'2:}(the Morse Smale condition) For any two critical manifolds
$\Sigma, \Sigma'$ and
$x\in \Sigma, \ $
$W_x^-$ and $W_{\Sigma'}^+$ are transversal.

As before
C'2 implies that
$\Cal M(\Sigma,\Sigma'):= W_{\Sigma}^- \cap W_{\Sigma'}^+$ is a
smooth
manifold of dimension equal to
 $i(\Sigma)- i(\Sigma') -\dim \Sigma,$ and
that the evaluation maps $u:\Cal M(\Sigma,\Sigma')\to \Sigma$ is a
smooth bundle
with fiber $W_x^-\cap W_{\Sigma'}^+$ a smooth manifold of dimension
$i(\Sigma) - i(\Sigma') .$

$\Cal M(\Sigma,\Sigma')$
is equipped with the free action
$\mu: \Bbb R\times \Cal M(\Sigma,\Sigma') \to
 \Cal M(\Sigma,\Sigma')$
defined by $\mu(t, z) = \gamma _z(t)$
and we denote the quotient space $\Cal M(\Sigma,\Sigma')/\Bbb R $
by $\Cal T(\Sigma,\Sigma').$
$ \Cal T(\Sigma,\Sigma')$ is a smooth manifold of dimension
$i(\Sigma)\ -
i(\Sigma') +\dim \Sigma -1,$
diffeomorphic to the submanifold
$ \Cal M(\Sigma,\Sigma') \cap h^{-1}(\lambda),$ for any real number
$\lambda$ in the open
interval $(h(\Sigma),  h(\Sigma')).$
In addition, one has the evaluation maps,
$u_{\Sigma, \Sigma'}:\Cal T(\Sigma,\Sigma')\to \Sigma$
which
is a smooth bundle
with fiber $W_x^-\cap W_{\Sigma'}^+/ \Bbb R,$ a smooth manifold of
dimension
$i(\Sigma) - i(\Sigma') -1,$ and
$l_{\Sigma,\Sigma'}:\Cal T(\Sigma,\Sigma')\to \Sigma'$
a smooth map. The maps $u_{\cdots}$ and $l_{\cdots}$ induce by pull-
back constructions the smooth bundles
$$u_{({\Sigma_0},{\Sigma}_1,\cdots, {\Sigma}_k)}:
\Cal T({\Sigma_0},{\Sigma}_1)\times
_{\Sigma_1}\cdots\times_{\Sigma_{k-1}}
\Cal T({\Sigma}_{k-1}, {\Sigma}_k)\to \Sigma_0,$$ the smooth maps
$$l_{({\Sigma_0},{\Sigma}_1,\cdots, {\Sigma}_k)}:
\Cal T({\Sigma_0},{\Sigma}_1)\times
_{\Sigma_1}\cdots\times_{\Sigma_{k-1}}
\Cal T({\Sigma}_{k-1}, {\Sigma}_k)\to \Sigma _k$$
and 
$$i_{\Sigma_0, \Sigma_1,\cdots ,\Sigma_k}: \Cal T(\Sigma_0,\Sigma_1)
\times_{\Sigma_1}
\cdots
\times_{\Sigma_{k-1}} \Cal T(\Sigma_{k-1},\Sigma_k)\times_{\Sigma_k} W_{\Sigma_k}^-
\to M,$$
defined by
$i_{\Sigma_0, y_1,\cdots ,y_k}(\gamma_1,\cdots, \gamma_k, y):=
i_{\Sigma_k}(y)$, for
$\gamma_i \in \Cal T(\Sigma_{i-1}, \Sigma_i)$ and
$y\in W_{\Sigma_k}^-,$

The analogue of Theorem 1.1 is Theorem 1.1' below.
The proof of Theorem 1.1 as given below, subsection d), and is
formulated in such way that the
extension to the Morse Bott case is straightforward.

\proclaim {Theorem 1.1'}

Let $\tau =(h,g)$  be a MB generalized triangulation.

1) For any critical manifold $\Sigma \subset Cr(h),$
the
smooth manifold $W^-_{\Sigma}$ has a canonical compactification
to a compact manifold
with corners $\hat W_{\Sigma}^-$, and the smooth bundle
$\pi^-_{\Sigma}:W^-_{\Sigma}\to \Sigma$
resp. the smooth inclusion
$i_{\Sigma}: W_{\Sigma}^- \to M$ have  extensions
$\hat{\pi}^-_{\Sigma}:\hat{W}^-_{\Sigma}\to \Sigma,$
a smooth bundle whose fibers are compact manifolds with corners,
resp. $\hat{i}_{\Sigma}: \hat{W}_{\Sigma}^- \to M$ a smooth map, so
that

\noindent (a):
$(\hat W_{\Sigma}^-)_k = \bigsqcup_{({\Sigma},{\Sigma}_1,\cdots,
 {\Sigma}_k)}
\Cal T({\Sigma},{\Sigma}_1)\times _{\Sigma_1}\cdots\times_{\Sigma_{k-1}}
\Cal T({\Sigma}_{k-1}, {\Sigma}_k)\times _{\Sigma_k}W_{{\Sigma}_k}^-,$

\noindent (b): the restriction of $\hat i_{\Sigma}$ to
$\Cal T({\Sigma},{\Sigma}_1)\times _{\Sigma_1}\cdots\times_{\Sigma_{k-1}}
\Cal T({\Sigma}_{k-1}, {\Sigma}_k)\times _{\Sigma_k}W_{{\Sigma}_k}^-$
is
given by $i_{\Sigma, \Sigma_1 \cdots, \Sigma_k}.$

2) For any two critical manifolds ${\Sigma}, {\Sigma'}$ with
$i(\Sigma)> i(\Sigma')$ the
smooth
manifold $\Cal T({\Sigma},{\Sigma'})$ has a canonical compactification
to a compact manifold
with corners $\hat {\Cal T}(\Sigma, \Sigma')$  and the smooth maps
$u:\Cal T(\Sigma, \Sigma') \to \Sigma$ and
$l:\Cal T(\Sigma,\Sigma')\to \Sigma'$ have smooth extensions
$\hat{u}:\hat{\Cal T}(\Sigma, \Sigma') \to \Sigma$ and
$\hat{l}:\hat{\Cal T}(\Sigma,\Sigma')\to \Sigma'$
with
$\hat {u}$ a smooth bundle whose fibers are compact manifolds with corners.
Precisely

$(\hat {\Cal T}(\Sigma,\Sigma'))_k =
 \bigsqcup_{({\Sigma},{\Sigma}_1,\cdots,
 {\Sigma}_k)}
\Cal T({\Sigma},{\Sigma}_1)\times
 _{\Sigma_1}\cdots\times_{\Sigma_{k-1}}
\Cal T({\Sigma}_{k-1}, {\Sigma}_k).$
\endproclaim

For a critical manifold $\Sigma$
choose an orientation of $\nu_-$ if this bundle is orientable and an
orientation of the orientable double cover of $\nu_-$ 
if not. Such an object will be denoted by
$\Cal O_\Sigma$ and the
collection of all $\Cal O_\Sigma$ will be denoted by
$\ o\equiv \{ \Cal O_\Sigma | \Sigma\subset Cr(h) \}.$

Choosing the collection $o$ in addition to  the Morse-Bott generalized
triangulation
$\tau$ one can provide as an analogue  of the geometric complex
$(C^*(M,\tau),\partial_{(\tau, o)}^*)$
the
complex $(C^\ast, D^\ast)$ defined by
$$C^r=: \bigoplus_{\{(k,\Sigma) | k+i(\Sigma)=r\}} \Omega ^k(\Sigma,
o(\nu_-))$$
and $D^r : C^r\to C^{r+1}$ given by the matrix
$||| \partial ^k_{\Sigma, \Sigma'}|||$ whose entries

$\partial ^k({\Sigma, \Sigma'}): \Omega^k(\Sigma',o(\nu'_-))\to \Omega
 ^{r-i(\Sigma)+i(\Sigma')+1}(\Sigma, o(\nu_-))$ are given by
\footnote{ this formula is implicit in (2.2) in view of the fact that 
$Int^\ast: (\Omega^\ast,d^\ast) \to (C^\ast(M,\tau),\partial^\ast)$
is supposed to be a surjective morphism of cochain complexes}

$$
\CD
\partial ^k({\Sigma, \Sigma'})= \left \{\aligned
&d^k:\Omega^k(\Sigma, o(\nu_-))\to \Omega
 ^{k+1}(\Sigma, o(\nu_-)) \text{if}\ \Sigma= \Sigma'\\
 &(-1)^k (\hat{u}_{\Sigma, \Sigma'})_{\ast}\cdot (\hat{
 l}_{\Sigma,
 \Sigma'})^{\ast}  \text {otherwise}
 \endaligned \right \}
 \endCD
 $$

Here
$\Omega^\ast(\Sigma, o(\nu_-))$ denotes the differential forms on
$\Sigma$ with
coefficients in the orientation bundle of $\nu_-$
and $(...)_{\ast}$ denotes the integration along the fiber of
$\hat{ u}_{\Sigma,\Sigma'}.$

The orientation bundle of $\nu_-$ has a canonical flat connection.
When $\nu_-$ is orientable then this bundle is trivial as bundle with
connection and $\Omega^\ast$ identifies to the ordinary differential
forms.

A Morse Bott function $h$ is a smooth function for which the critical
set consists of a
disjoint
union of connected manifolds $\Sigma,$ so that the Hessian of $h$ at
each critical point
of
$\Sigma$ is nondegenerated in the normal direction.

{\bf Observation:}

O.2': Given a pair $(h,g)$ with $h$ a Morse Bott function and $g$ 
is a Riemannian metric
one can provide arbitrary small $C^0$ perturbation $g'$ of $g$  so
that the pair
$(h, g')$
satisfies $C'.$ If $(h,g)$ satisfies $C'1
$ one can choose $g'$ arbitrary closed to $g$ in $C^0-$topology
so that $g=g'$ away from a given neighborhood
of the critical point set, $g'=g$ in some (smaller) neighborhood of the 
critical point set and 
$(h,g')$ is an MB generalized triangulation.

\vskip .2in

\subhead {\bf  c. G-Generalized triangulations}
\endsubhead

\vskip .2in
Of particular interest is the case of a smooth $G-$manifold
$(M, \mu:G\times M\to M)$
where $G$ is a compact Lie group and $\mu$ a smooth action.
In this case we consider pairs $(h, g)$ with $h$ a G-invariant smooth
function and $g$ a $G-$invariant Riemannian metric. Then
$Cr(h)$ consists of
a union of $G-$orbits. The $G-$version of conditions {\bf C1} and {\bf
C2}
are obvious to formulate.
We say that the pair $(h,g)$ with both $h$ and $g\ $ $G-
$invariant is a {\it G-generalized triangulation }resp.
{\it normal G-generalized triangulation} if G-C1 (resp. normal G-C1) and
G-C2 hold.

{\bf G-C1:} $Cr(h)$ is a finite union of orbits denoted by $\Sigma $  and for
any critical orbit
 $\Sigma$ we require the existence of an admissible chart. More precisely,
 such a chart
 around $\Sigma$ is provided
by the following data:

1: A closed subgroup $H\subset G,$  two orthogonal representations
$\rho_{\pm}:H \to O(V_{\pm})$ and a scalar product on the Lie algebra
$\frak g$
of $G$ which is invariant with respect to the adjoint representation
restricted to $H.$

$\rho_{\pm}$ induce orthogonal bundles
$\nu_{\pm}: E(\nu_{\pm})\to \Sigma.$
The total space of these $G-$ bundles
are $E(\nu_{\pm})= G\times _H V_{\pm}.$ The scalar product on
$\frak g$ and the
the scalar product on $V= V_-\oplus V_+$ induces a
$G-$invariant Riemannian metric on $G$ and on $G\times V.$ The metric
on
$G\times V$ descends to a $G-$invariant  Riemannian metric on
$E(\nu_-\oplus\nu_+).$

2) A positive number $\epsilon,$ a  constant $c\in \Bbb R$ and a  $G-
$equivariant diffeomorphism
$\varphi: (U,\Sigma) \to D_{\epsilon}(\rho_+\oplus \rho_-)$ where $U$
is a closed  $G-$
tubular neighborhood  of $\Sigma$ in $M,$ and $D_{\epsilon}$ denotes the 
disc of of radius $\epsilon$ in the underlying Euclidean space of the representation
$\rho_+\oplus \rho_-,$ so that

$(i): h\cdot \varphi^{-1}((g,v_1,v_2))=c - 1/2 ||v_1||^2+ 1/2
 ||v_2||^2) $
where $v_\pm\in E(\nu_\pm)$

$(ii): (\varphi^{-1})^*(g)$ is the Riemannian metric on
$E(\nu_-\oplus\nu_+)$ described above.

We call the admissible chart "normal" if in addition $\rho_-$ is
trivial.
The condition {\bf normal G-C1} requires the admissible charts to be
normal.

{\bf G-C2:} This condition is the same as C'2.

{\bf Observation: }

O1": Given a pair $(h,g)$  one can perform an arbitrary small $C^0$
perturbation $(h',g')$ 
so that $(h', g')$ satisfies normal G-C1. This was
proven in [M].

O2": Given $(h,g)$ a pair which satisfies normal G-C1
then one can perform an arbitrary small $C^0$ perturbation on
the metric $g$ and obtain the $G-
$ invariant metric $g'$
(away from the critical set) so that $(h,g')$ satisfies G-C2.
This result is proven in [B].

 Clearly, a $G-$generalized triangulation is a $MB-
 $generalized triangulation, hence Theorem 1.1'
 above can be restated in this case as Theorem 1.1" with the
 additional specifications
 that in the statement of Theorem 1.1' all compact
 manifolds with corners are $G-$ manifolds and all maps are $G-$
 equivariant.

Note that  a $G-
$generalized triangulation provides via Theorem 1.1" a structure of a
smooth $G-$handle
body and
a normal $G-$ generalized triangulation provides a structure of smooth 
$G-$CW complex for $M.$
The smooth triangulability of
compact smooth manifolds with corners if combined with the existence
of normal $G-
$ generalized
triangulation lead to the existence of a smooth $G-$triangulation in
the sense of [I] and
then to the existence of smooth
triangulation of the orbit spaces of a smooth $G-$manifold when $G$
is
compact. There is no proof for this result in literature. The best
known
result so far, is the existence of a $C^0$ triangulation of the orbit space
established
by Verona.
[V].

\vskip .2in

\subhead {\bf d): Proof of Theorem 1.1}
\endsubhead
\vskip .1in

{\it Some notations}

We  begin by introducing some notations:

\noindent Let  $c_0< c_1\cdots < c_N$ be the collection of all
critical values ($c_0$
the absolute minimum, $c_N$ the absolute maximum) and fix
$\epsilon >0$
small enough so that $c_i-\epsilon > c_{i-1}+ \epsilon$ for all
$i\geq 1.$
Denote by:
$$Cr(i):= Cr(h)\cap h^{-1}(c_i),$$
$$M_i:= h^{-1}(c_i),$$
$$M^{\pm}_i:=h^{-1}(c_i \pm \epsilon)$$
$$M(i):=h^{-1}(c_{i-1}, c_{i+1})$$
\vskip .1in

For any $x\in Cr(i)$ denote by:
$$S^{\pm}_x:= W_x^{\pm}\cap\   M^{\pm}_i$$
$${\bold S}_x:= S^+_x \times S^{-}_x$$
$$W_x^{\pm}(i):= W_x^{\pm}\cap M(i)$$
$${\bold SW}_x(i):= S_x^+\times W_x^-(i).$$

\noindent It will be convenient to write
$$S^{\pm}_i:= \bigcup _{x\in Cr(i)}S^{\pm}_x$$
$$\bold S_i:= \bigcup _{x\in Cr(i)}S_x \subset M^{-}_i
\times
M^{+}_i$$
$$W^{\pm}(i):= \bigcup _{x\in Cr(i)}W_x^{\pm}(i)$$
$$\bold SW(i):= \bigcup _{x\in Cr(i)}S^{+}_x\times W^{-}_x(i)$$

\dofig{4in}{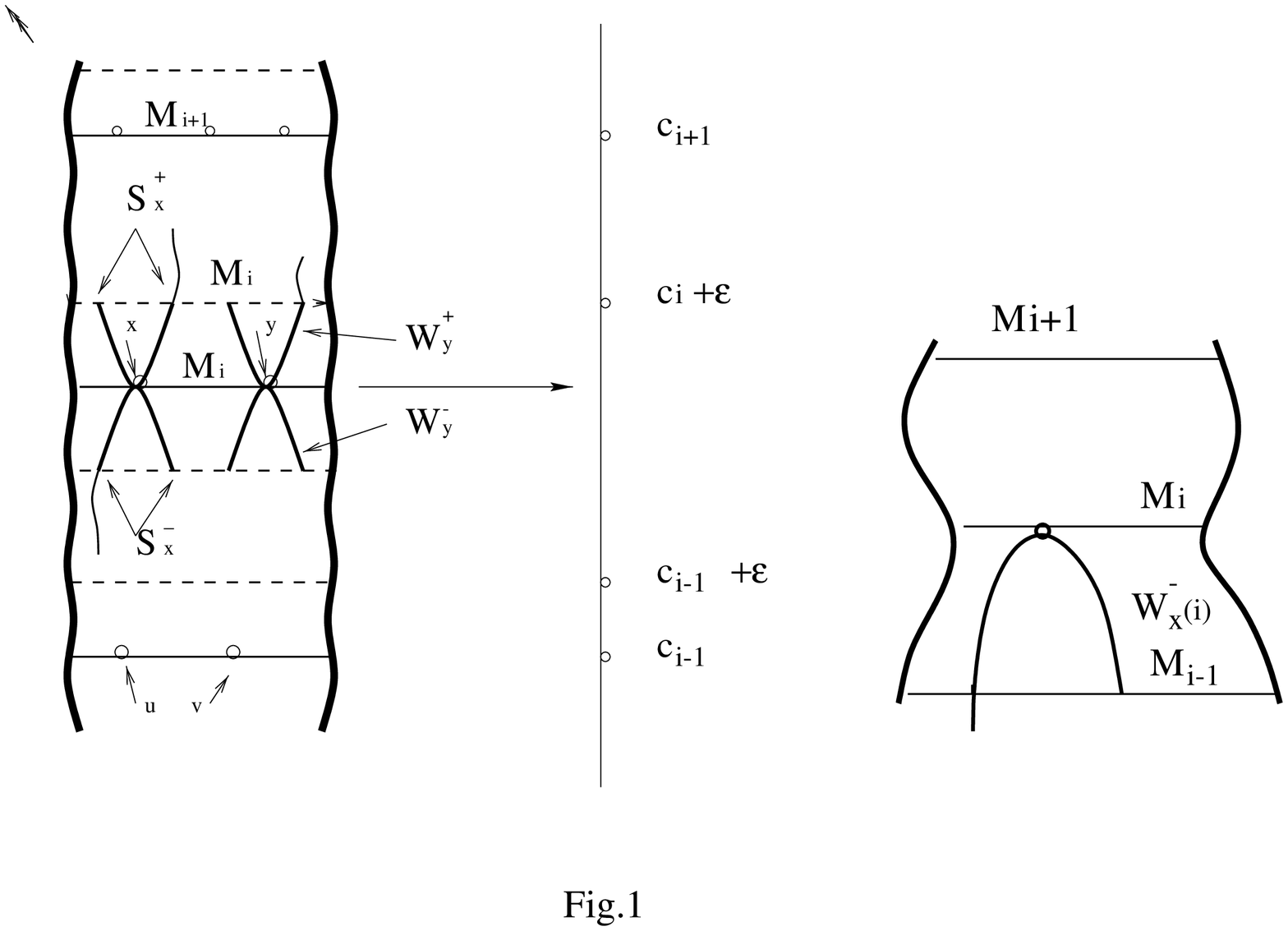}

\noindent Observe that:

1)
$\bold S_i \subset M^+_i \times M^-_i,\ \ \bold SW(i)\subset
 M^+_i\times M(i)$

2) $ M^{\pm}_i$ is a smooth
manifold of
dimension $n-1$, ($n=\dim M$) and
$M(i)$ is a smooth manifold of dimension $n,$ actually an open
set in $M.$
$M_i$ is not a
manifold, however,
$\overset\circ \to M_i: =M_i\setminus Cr(i), \ \
{\overset\circ \to M}_i^{\pm}: = M_i^{\pm}\setminus
S_i^{\pm}$
are
are smooth  manifolds (submanifolds of $M$) of dimension $n -1.$
\vskip .2in

{\it The flow $\Phi_t$ and few induced maps}

Let  $\Phi_t$ be the flow associated to the vector field
$-\text{grad}_g h/||-\text{grad}_g h||$ on

\noindent $M\setminus Cr(h)$
and consider:

a) the diffeomorphisms
$$\psi_i: M^-_{i}\to M^+_{i-1}$$
$$\varphi_i^{\pm}: {\overset\circ \to {M^{\pm}} _i}\to {
\overset\circ \to M}_i$$
obtained by the restriction of
$\Phi_{(c_i -c_{i-1} -2\epsilon)}$ and $\Phi_{\pm\epsilon} ,$

b) the submersion
$\varphi(i):M(i)\setminus (W^-(i)\cup W^+(i)) \to \overset\circ \to
 M_i$
defined by $\varphi(i) (x):= \Phi_{h(x)-c_i}(x).$

Observe that $\varphi_i^{\pm}$ and $\varphi(i)$ extend in an unique way to 
continuous maps
$$\varphi_i^{\pm}: M^{\pm}_i\to M, \ \ \varphi(i):M(i)\to M_i.$$

\dofig{4in}{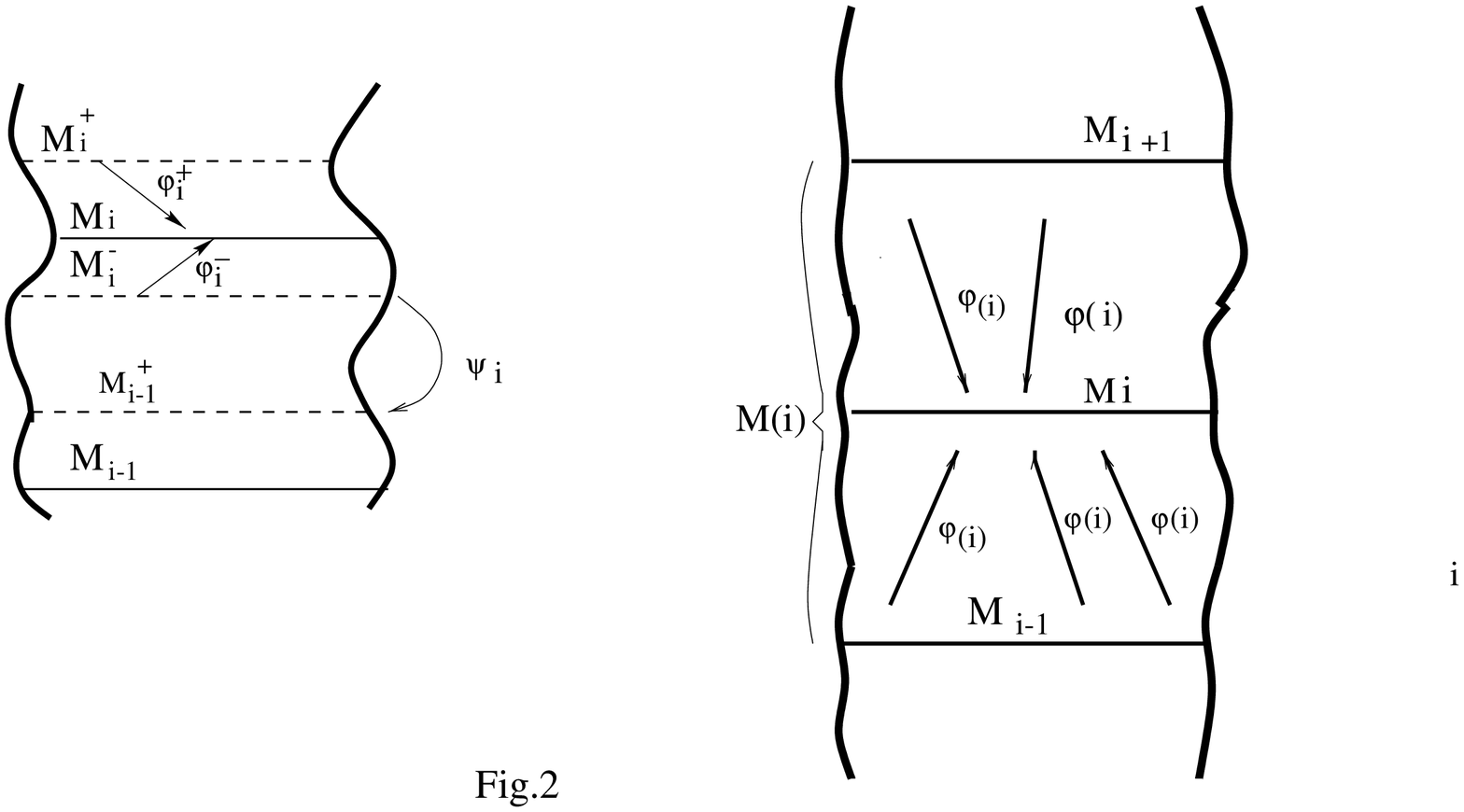}

\vskip .2in
{\it Two manifolds with boundary}
\vskip .1in
{\it The manifold $P_i:$}
Define
$$P _i:=\{(x,y)\in  M_i^{-}\times  M_i^{+}
|\ \varphi^-_i(x)=
\varphi^+_i(y)\},$$
and denote  by $p^{\pm}_i:
P_i\to M_i^{\pm}$ the canonical projections.
One can verify that $P_i$ is a compact smooth $(n-1)$ dimensional
manifold with boundary, (smooth submanifold of $M^-_i\times M^+_i$)
whose boundary $\partial P_i$ can be identified to
$\bold S_i\subset M^-_i\times M^+_i.$ Precisely

OP1:
$p^{\pm}_i: P_i\setminus \partial P_i \to
{\overset \circ \to {M}} _i^{\pm}$ are diffeomorphisms,

OP2: the restriction of $p^+_i\times p^-_i$ to $\partial P_i$
is a diffeomorphism onto $\bold S_i.$ (Each $p^{\pm}_i$ restricted
to $\partial P_i$ is the projection onto
$S^{\pm}_i.$)
\vskip .1in

{\it The manifold $Q(i):$} Define
$$Q(i)=\{(x,y)\in  M_{i}^{+}\times  M(i) |
\varphi^+_i(x)= \varphi(i)(y)\}$$
or equivalently, $Q(i)$
consists  of pairs of
points $(x, y),$ $x\in M_{i}^{+}$, $y\in M(i)$
which lie on the
same
(possibly broken) trajectory
and denote  by $l_i:Q(i)\to M^+_i$ resp.
$r_i:Q(i)\to M(i)$  the canonical projections.

One can verify that $Q(i)$ is a  smooth $n-$dimensional
manifold with boundary, (smooth submanifold of $M^+_i\times M(i)$)
whose boundary $\partial Q(i)$ is diffeomorphic to
$\bold SW(i)\subset M^+_i\times M(i).$ More precisely

OQ1:
$l_i: Q(i)\setminus \partial Q(i) \to
{\overset \circ \to {M}}^+ _i$ is a smooth bundle with
fiber an open segment
and $r_i: Q(i)\setminus \partial Q(i) \to M(i)\setminus W^-(i)$
a  diffeomorphism,

OQ2: the restriction of $l\times r$ to $\partial Q(i)$
is a diffeomorphism onto $\bold SW(i).$ ($l$ resp. $r$ restricted
to $\partial Q(i)$ identifies with the projection onto
$S^{+}_i$ resp. $W(i)$).
\vskip .1in

Since $P_i$ and $Q(i)$ are smooth manifolds with
boundaries, 
$$\Cal P_{r, r-k}:=P_r\times P_{r-1}\cdots P_{r-k}$$ and
$$\Cal P_r(r-k):=P_r\times \cdots P_{r-k+1}\times Q(r-k)$$
are smooth manifolds with corners.

Our arguments for the proof of Theorem 1.1 will be based on the
following 
method for recognizing a smooth manifold with corners.
If
$\Cal P$ is a smooth manifold with corners, $\Cal O, \Cal S$ smooth
manifolds,
$p:\Cal P \to \Cal O$ and $s:\Cal S \to \Cal O$ smooth maps so that
$p$ and
$s$ are transversal
($p$ is transversal to $s$ if its restriction to each $k-$ corner of
$\Cal P$ is
transversal to $s$, ) then $p^{-1}(s(\Cal S))$ is a smooth submanifold
with corners of $\Cal P.$

{\bf Proof of Theorem 1.1 } First we prove part (2). We want to verify
that $\hat{\Cal T}(x,y)$ (cf the definition in the statement of
Theorem 1.1) is a smooth manifold with corners. Let
$x\in Cr(r+1)$ and
$y\in Cr(r-k-1), \ k\geq -2$.
If $k=-2$ the statement is empty, if $k=-1$ there is nothing
to check, so we suppose $k\geq 0.$

We consider
$\Cal P= \Cal P_{r,k}$ as defined above,
$\Cal O=\prod^{r-k}_r (M_i^+\times M_i^-))$, and
$\Cal S = S^-_x\times M^-_r\cdots \times M^-_{(r-k+1)}\times S^+_y.$
In order to define the maps $p$ and $s$ we consider

$\omega_i: M^-_i \to
 M^{-}_i\times  M^{+}_{i+1}\ $ given by
$\ \omega_i(x)=(x, \psi_i(x))$,
and

$\tilde p_i:P_i\to  M^{+}_i\times  M^{-}_i\ $
given by
$\ \tilde p_i(y)=
(p^+_i(y), p^-_i(y))$.

We also denote  by $\alpha: S^-_x\to  M^+_r $
and $\beta: S^+_y\to  M^-_{k-r}$ the restriction of
$\psi_{r+1}$
resp. of $\psi_{r-k}^{-1}$ to $S_x^-$ resp. $S_y^+.$
Take $s= \alpha\times \omega_r\cdots \omega _{r-k+1}\times \beta$
and $p:=\prod_{i=r}^{r-k} \tilde \pi_i.$

The verification of the transversality follows easily from
OP1, OP2 and the Morse Smale
condition C2.
It is easy to see that $p^{-1}(s(\Cal S))
$ is compact and identifies to
$\hat {\Cal T}(x,y);$ the verification of this fact is
left to the reader.

To prove part (1) we first consider the map $\hat i_x: X=\hat W^-_x \to M$
defined by (a) and (b) in Theorem 1.1 (1). Let $X:=\hat W^-_x$ and  
for any positive integer $k$ we denote by $X(k):= \hat i_x^{-1}(M(k)).$
The proof will be given in two steps. First we will topologize
$X(k)$ and put on it
a structure of a smooth manifold with corners,
so that the restriction of $\hat i_x$ to $X(k)$ is a smooth map.
Second we check that $X(k)$ and $X(k')$ induce on
the intersection $X(k)\cap X(k')$ the
same topology and the same smooth structure.
These facts imply that $X$ has a canonical structure of
smooth manifold with corners and $\hat i_x$ is a smooth
map. The compacity of $X$ follows  by observing that
the image $\hat i_x(X)$ is compact and the preimage of any point is
compact.

To accomplish the first step we proceed in exactly the same way as in the
proof of part (2).
Suppose $x\in Cr(r+1).$ Consider

$\Cal P:= \Cal P_r (k),$

$\Cal O:=(M_r^+\times
 M_{r}^-)\times \cdots ( M_{r-k+1}^+\times
 M_{k+r-1}^-)\times  M_{r-k}^+$, and

$\Cal S:= S^-_x\times M^-_r\cdots \times M^-_{r-k+1}.$

\noindent Take

$p:=\tilde p_r\times \cdots \times \tilde p_{r-k+1}\times
l_{r-k}$
and

$s= \alpha\times \omega_r\cdots \times \omega_{r-k+1}.$

The verification of the transversality follows from OP1, OP2, OQ1, OQ2
and the Morse
Smale condition C2
above as explained in the Appendix. It is easy to see that
$p^{-1}(s(\Cal S))$ identifies to
$X(r-k).$

The second step is more or less straightforward, so it will be left
again to the reader.
q.e.d.

\subhead {Appendix : The verification of transversality of the 
maps $p$ and $s$}
\endsubhead
\vskip .1in
Consider the diagrams

$$
\xymatrix @C=3pt
{S^-_{r+1}\ar[d]^{\alpha} & & &M^-_r\ar[dl]^{id}\ar[dr]^{\psi_r} & & &
&{M^-_{r-1}}\ar[dl]^{id} &
\cdots &{M^-_{r-k+1}}\ar[dr]^{\psi_{r-k+1}} & & &S^+_{r-k-1}\ar[d]^{
\beta}\\
M^+_r\ar[dr] & &{M^-_{r}} & &{M^+_{r-1}}& &M^-_{r-1} & &\cdots &  &{M^
+_{r-k}} & &{M^-_{r-k}} \\
 & P_r\ar[ul]^-{p^+_{r}}\ar[ur]^{p^-_{r}} & & & &P_{r-1}\ar[ul]^{p^+_{
 r-1}}\ar[ur]^{p^-_{r-1}} & & &
 \cdots & & &P_{r-k}\ar[ul]^{p^+_{r-k}}\ar[ur]^{p^-_{r-k}}}
$$
\hskip 2.5in Diagram 1
\vskip .2in

$$
\xymatrix @C=3pt
{S^-_{r+1}\ar[d]^{\alpha} & & &{M^-_r}\ar[dl]^{id}\ar[dr]^{\psi_r} & &
&
&{M^-_{r-1}}\ar[dl]^{\psi_{r-k+1}} &\cdots &{M^-_{r-k+
1}}\ar[dr]^{\psi_{r-k+1}} & \\
{M^+_r} & &{M^-_{r}} & &{M^+_{r-1}}& &M^-_{r-1} & &\cdots &  &{M^+_{r-
k}}\\
 &{P_r}\ar[ul]^{p^+_r}\ar[ur]^{p^-_r} & & &
 &P_{r-1}\ar[ul]^{p^+_{r-1}}\ar[ur]^{p^-_{r-1}} & &
 &\cdots & &Q(r-k)\ar[u]_{l_{r-k}}}
$$
\hskip 2.2in Diagram 1'
\vskip .2in

$$
\xymatrix @C=3pt
{ &M^-_r\ar[dl]^{id}\ar[dr]^{\psi_r} & & & &{M^-_{r-1}}\ar[dl]^{id} &
\cdots
&{M^-_{r-k+
1}}\ar[dr]^{\psi_{r-k+1}} & & &S^+_{r-k-1}\ar[d]^{\beta}\\
{M^-_{r}} & &{M^+_{r-1}}& &M^-_{r-1} & &\cdots &  &{M^+_{r-k}} & &{M^-
_{r-k}} \\
S^-_{r}\ar[u]^i & & &P_{r-1}\ar[ul]^{p^-_{r-1}}\ar[ur]^{p^+_{r-1}} & &
&\cdots & &
&P_{r-k}\ar[ul]^{p^+_{r-k}}\ar[ur]^{p^-_{r-k}}}
$$
\hskip 2.2in Diagram 2
\vskip .2in

$$
\xymatrix @C=3pt
{S^-_{r+1}\ar[d]^\alpha & & &{M^-_r}\ar[dl]^{id}\ar[dr]^{\psi_r} & & &
&{M^-_{r-1}}\ar[dl]^{id}
&\cdots &{M^-_{r-k+
1}}\ar[dr]^{\psi_{r-k+1}} & \\
{M^+_r} & &{M^-_{r}} & &{M^+_{r-1}}& &M^-_{r-1} & &\cdots &  &{
M^+_{r-k}} \\
 &{P_r}\ar[ul]^{p^+_r}\ar[ur]^{p^-_r} & & &
 &P_{r-1}\ar[ul]^{p^+_{r-1}}\ar[ur]^{p^-_{r-1}} & & &\cdots
 & &S^+_{r-k}\ar[u]^i}
$$
\hskip 2.2in Diagram 3

\vskip .2in

$$
\xymatrix @C=3pt
{ &{M^-_r}\ar[dl]^{id}\ar[dr]^{\psi_r} & & & &{M^-_{r-1}}\ar[dl]^{id}
&\cdots &{M^-_{r-k+
1}}\ar[dr]^{\psi_{r-k+1}} & \\
{M^-_{r}} & &{M^+_{r-1}}& &M^-_{r-1} & &\cdots &  &{M^+_{r-k}}
\\
S^-_{r}\ar[u]^i & & &P_{r-1}\ar[ul]^{p^+_{r-1}}\ar[ur]^{p^-_{r-1}} & &
&\cdots & &Q(r-k)\ar[u]
^{l_{r-k}}}
$$
\hskip 2.2in Diagram 4
\vskip .2in

For each of these diagrams denote by $\Cal P$ resp. $\Cal O$ resp.
$\Cal R,$ the product of the manifolds on the third  resp. second
resp. first row and let $p: \Cal P \to \Cal O$ resp.
$s:\Cal R \to \Cal O$ denote the product of the maps from the third to
the second row resp. the first to the second row.
Clearly $\Cal P$ is a smooth manifold with corners. Denote by
${\overset \circ \to {\Cal P}}$
the interior of $\Cal P,$ i.e $\Cal P \setminus  \partial \Cal P,$
and by
${\overset \circ \to {p}}: {\overset \circ
\to {\Cal P}}\to \Cal O$ the restriction of $p$
to ${\overset \circ \to {\Cal P}} .$

We refer to the  statement " ${\overset \circ \to {p}}$
transversal to $s$" with
$p$ and $s$ obtained from the diagram 1, 1' 2, 2',
3, 4 as:
$T^1_{r, k},\  T^{1'}_{r, k},\  T^2_{r, k},\  T^{2'}_{r, k},\
T^3_{r, k},\  T^4_{r, k}.$

In view of the observation OP2 and OQ2
and of the fact that "
$m_i:\Cal M_i\to \Cal O_i \text{and} s_i:\Cal R_i\to \Cal O_i,$ $i=1,2,$
$\Cal M_i, \Cal O_i, \Cal R_i$
smooth manifolds imply the transversality of
$m_1\times m_2:\Cal M_1\times \Cal M_2\to \Cal O_1
\times \Cal O_2 \text{and}\  s_1\times s_2: \Cal R_1
\times \Cal R_2\to \Cal O_1\times \Cal O_2$",
it is easy to see that the transversality of $p$ and  $s$ obtained
from the diagram $1$ resp. $1'$ can be derived
from the validity of
the statements
$T^1_{r, k},\  T^2_{r, k},\  T^3_{r, k}$ resp.
$T^{1'}_{r, k},\
T^3_{r, k},\  T^4_{r, k}$ for various $r, k.$

In view of the fact that all arrows except
$\alpha, \beta, l_{r-k}$  and the inclusions ( cf Diagrams 2,3,4)  are
open embeddings,  the properties
$T^1_{r, k},\ \cdots,  T^{4}_{r, k}$ follow
from
the transversality of $W^-_{r+1}$ and $W^+_{r-k-1}.$
This finishes the verification of the transversality
statement needed for the proof of Theorem 1.1. \hskip 1in
q.e.d
\vskip .5in

\proclaim{ Lecture 2: Witten deformation and the spectral properties of the
Witten
Laplacian }\endproclaim

\vskip .2in

\subhead {a. De Rham theory and integration}
\endsubhead

 Let $M$ be a closed smooth  manifold and $\tau= (h,g) $ be a
 generalized
 triangulation or more general a Morse-Bott  generalized triangulation
 and $ o$ a
 system of orientations.
 Denote by $(\Omega^*(M), d^*)$ the De Rham complex of $M.$ This is a
cochain complex whose component $\Omega^q(M)$ is the (Frechet) space
of
 smooth differential forms of degree $q$ and whose differential
$d^q: \Omega^q(M) \to \Omega^{q+1}(M)$
is given by the exterior differential $d.$ Recall that Stokes' theorem for manifolds with corners
can be formulated
as follows:

 \proclaim{Theorem } Let $P$ be a compact $r-$dimensional oriented
smooth manifold with corners
 and $f:P\to M$  a smooth map. Denote by $\partial f:
P_1\to M $ the restriction of $f$ to
 the smooth oriented manifold $P_1$ ($P_1$ defined as in section 1a).
If $\omega \in \Omega^{r-1}(M)$
 is a smooth form then $\int_{P_1} (\partial f)^*(\omega) $
is convergent and
 $$\int_P f^*(d\omega)= \int_{P_1} (\partial f)^*(\omega).$$
 \endproclaim

Define the linear maps
$Int^q: \Omega^q(M) \to C^q(M,\tau)$ as follows.

a) In the case $\tau$ is a generalized triangulation and
$C^q(M,\tau):=
\text{Maps}(Cr_q(h),\Bbb K),$
$Int^q$ is defined by
$$ \leqalignno{Int^q(\omega)(x):= \int _{\hat W_x^-} \omega,&&(2.1)\cr}
$$
where $\Bbb K$ is the field of real or complex numbers.

b) In the case $\tau$ is an Morse-Bott generalized triangulation and 
$$C^q(M,\tau):= \bigoplus_{k+\dim(\Sigma)=q} \Omega ^k
(\Sigma,o(\nu_-)),$$
$Int^q$ is given by integrating the pull back by $\hat{i}_{\Sigma} :\hat{W}^-_{\Sigma} \to M$
of a form $\omega \in \Omega^q(M)$ along the fiber of 
$\hat{\pi}_\Sigma :\hat{W}^-_\Sigma \to \Sigma,$  i.e.
$$\leqalignno{Int^q(\omega)= \oplus_{k+\dim(\Sigma)=q}(\hat{\pi}^-_{\Sigma})_{\ast} \hat{i}_
{\Sigma}^\ast \omega \in
  \Omega^{q-i(\Sigma)}(\Sigma,o(\nu_-))&&(2.2)\cr}$$

It is a consequence of Theorem 1.1 (1.1') that
the collection of the linear maps ${Int^q}$ defines a morphism
$$Int^*: (\Omega^*(M), d^*) \to (C^*(M,\tau),\partial^*)$$
of cochain complexes which has the following property (de Rham):

\proclaim {Theorem }$Int^*$  induces an
isomorphism in cohomology.
\endproclaim
\vskip .1in

c) Suppose that we are in the case of a $G-$manifold and $\tau$ is a $G-
$generalized triangulation
as described in section 1.
Fix an irreducible representation $\xi: G \to O(V_\xi)$ and denote by
$(\Omega^\ast_\xi, d^\ast_\xi)$ resp. $(C^\ast_\xi, \partial^\ast_\xi)$ the
subcomplex of
$(\Omega^\ast, d^\ast)$ resp.
$(C^\ast, \partial^\ast)$ defined by the property that $\Omega^\ast_\xi$
resp. $C^\ast_\xi
$ is the largest $G-$invariant subspace
of
$\Omega^\ast$ resp. $C^\ast$ which contains no other irreducible representation but 
$\xi.$ Equivalently for any irreducible
representation $\xi',\ \xi'\ne
\xi,$
$\Omega^\ast_\xi \otimes_G V_{\xi'}=0\ $resp.
 $C^\ast_\xi\otimes_G V_{\xi'}= 0.$
Since $G$ is compact, 
$(\Omega^\ast, d^\ast)=\hat{\oplus}_\xi (\Omega^\ast_\xi, d^\ast_\xi)$
resp. $(C^\ast,
\partial^\ast)=
\hat{\oplus}_\xi (C^\ast_\xi , \partial^\ast_\xi).$
Since the integration map is $G-$equivariant, 
$Int(\Omega^\ast_\xi) \subset C^\ast_\xi.$

Denote by
$$Int^\ast_\xi: (\Omega^\ast_\xi, d^\ast_\xi)\to (C^\ast_\xi, 
\partial^\ast_\xi)$$
the restriction of $Int^\ast$ to the
components $(\Omega^\ast_\xi,d^\ast_{\xi})$ corresponding to $\xi.$

We have the following refinement of de Rham Theorem.

\proclaim {Theorem }$Int^*_\xi$  induces an
isomorphism in cohomology.
\endproclaim

Theorem 3.1-1" below (Lecture 3) will imply these theorems.

\vskip .2in

\subhead {b) Witten deformation and Witten Laplacians}
\endsubhead

\vskip .2in

Let $M$ be a closed manifold and $\alpha \in \Omega^1(M)$
a closed
$1-$ form, i.e $d\alpha=0.$
For $t>0$
consider
the complex $(\Omega^*(M), d^*(t))$
with differential
$$
\leqalignno{ d^q(t)(\omega)= d\omega +t \alpha \wedge \omega.
&&(2.3)\cr}
$$

If $\alpha = dh$ with $h:M\to \Bbb R$ a smooth function, 
$d^q(t)(\omega) = e^{-th}d e^{th}(\omega)$ and
$d^*(t)$ is the unique differential in $\Omega^*(M)$ which makes the
multiplication
by the smooth function $e^{th}$
an isomorphism of cochain complexes
$$e^{th}: (\Omega^*(M), d^*(t)) \to (\Omega^*(M), d^*).$$

Recall that for any vector field $X$ on $M$ one defines a zero
order differential operator,
$\iota_X^*: \Omega^*(M) \to \Omega^{*-1}(M),$ by
$$\leqalignno{(\iota_X^q\omega) (X_1, X_2,\cdots, X_{q-1}) := \omega (X,
X_1,\cdots, X_
{q-1})
&&(2.4)\cr}$$
and a first order differential operator
$L_X^*:\Omega^*(M) \to \Omega^*(M),$ the Lie derivative in the
direction $X,$ by
$$\leqalignno{ L_X^q:= d^{q-1}\cdot \iota_X^q + \iota_X^{q+1}\cdot d^
q.
&&(2.5)\cr}$$
They satisfy the following identities:
$$\leqalignno{\iota_X(\omega_1\wedge \omega_2)= \iota_X(\omega_1)
\wedge \omega_2 +
(-1)^{|\omega_1|}\omega_1 \wedge \iota_X(\omega_2)
.&&(2.6)\cr}$$
where $|\omega_1|$ denotes the degree of $\omega_1$ and
$$
\leqalignno{ L_X(\omega_1\wedge \omega_2)= L_X(\omega_1)\wedge \omega_
2 +
\omega_1 \wedge L_X(\omega_2).
&&(2.7)\cr}$$

Given a Riemannian metric $g$ on the oriented manifold $M$ we have the
zero
order operator
$R^q:\Omega^q(M) \to \Omega^{n-q}(M)$, known as the Hodge-star
operator
which, with respect to an oriented orthonormal frame
$e_1, e_2,\cdots, e_n$ in the cotangent space at $x,$ is given by
$$\leqalignno{ R_x^q (e_{i_1}\wedge\cdots \wedge e_{i_q}) = \epsilon (
i_1,\cdots,i_q)
e_1\wedge \cdots\wedge
\hat e_{i_1}\wedge\cdots \wedge\hat e_{i_q}\wedge\cdots\wedge e_n
,&&(2.8)\cr}$$
where $ 1\leq i_1<i_2 <\cdots <i_q\leq n ,$ and
$\epsilon (i_1,i_2,\cdots,i_q)$
denotes the sign of the
permutation of $(1,\cdots, n)$ given by
$$(i_1,\cdots,i_q, 1,2,\cdots, \hat i_1,\cdots,\hat i_2,\cdots, \hat
i_2,\cdots,
\hat i_q,\cdots,n).$$
Here a ``hat'' above a symbol means the deletion of this symbol.

The operators $R^q$ satisfy
$$
\leqalignno{ R^q\cdot R^{n-q}= (-1)^{q(n-q)}Id.&&(2.9)\cr}$$

With the help of the operators $R^q$ for an oriented Riemannian
manifold
of dimension $n,$ one defines the fiberwise scalar product
 $\ll\ ,\ \gg:\Omega^q(M)\times \Omega^q(M) \to \Omega^0(M)$ and the
 formal adjoints
of
$d^q, d^q(t), \iota_X^q, L_X^q,$
by the formulas
$$\leqalignno{ \ll \omega_1, \omega_2 \gg= (R^n)^{-1}(\omega_1\wedge R
^q
(\omega_2)),&&(2.10)\cr}$$
$$\leqalignno{
\delta^{q+1} & = (-1)^{nq+1} R^{n-q}\cdot d^{n-q-1}\cdot R^{q+1}:
\Omega^{q+1}(M) \to \Omega^q(M),\  &(
2.11)\cr
\delta^{q+1} &(t) = (-1)^{nq+1} R^{n-q}\cdot d^{n-q-1}(t)\cdot
R^{q+1}:\Omega^{q+1}(M) \to \Omega^q(M), \cr
(\iota_X^q)^{\sharp}&= (-1)^{nq-1} R^{n-q}\cdot \iota_X^{n-q-1}\cdot R
^{q-1}:\Omega^{q-1}(M) \to \Omega^q(M),\cr
(L_X^q)^{\sharp}&= (-1)^{(n+1)q+1} R^{n-q}\cdot L_X^{n-q}\cdot R^{q}
:\Omega^{q}(M) \to \Omega^q(M)\cr}
$$
These operators  satisfy
$$\leqalignno{\ll \iota_X \omega_1, \omega_2\gg\   =\  \ll \omega_1, (
\iota_X)^{\sharp}
\omega_2 \gg &&(2.12)\cr}$$
and
$$\leqalignno{ (L_X)^{\sharp} = (\iota_X )^{\sharp}\cdot \delta +
\delta \cdot
(\iota_X)^{\sharp}.
&&(2.13)\cr}$$

Note that
$ L_X^q + (L_X^q)^{\sharp}$ is a zeroth order differential operator.
Let $X^{\sharp}$ denote the element in $\Omega^1(M)$ defined by
$X^{\sharp}(Y):= \ll X, Y \gg\ $
and for  $\alpha \in \Omega^1(M)$
let $(E_{\alpha})^q: \Omega^q(M) \to \Omega^{q+1}(M),$ denote
the exterior product by $\alpha.$ Then we have
$$\leqalignno{(\iota_X^q)^{\sharp}= \ (E_{X^{\sharp}})^{q-1}.
&&(2.14)\cr}$$

It is easy to see that the scalar products $\ll , \gg$ and hence the
operators
$\delta^q ,\delta^q(t), \iota_X^{\sharp}$ and $ L_X^{\sharp}$ are
independent
of
the orientation of $M.$  Therefore they are defined (first locally and
then, being
differential operators, globally) for an arbitrary
Riemannian
manifold, not necessarily
orientable, and satisfy
 (2.10), (2.12)-(2.14) above.

For a Riemannian manifold $(M,g)$ one introduces the scalar product
$\Omega^q(M)\times \Omega^q(M)\to \Bbb C$
by
$$\leqalignno{< \omega, \omega'> := \int_M \omega \wedge \omega' =
\int_M \ll \omega, \omega'\gg dvol(g). &&(2.15)\cr}$$

In view of (2.12), $\delta^{q+1}(t), (\iota_X^q)^\sharp$ and
$(L_X^q)^\sharp$
are
formal adjoints of  $d^q(t), \ \iota_X^q(t)$ and $L_X^q$ with respect
to the
scalar product $<.,.>.$

For a Riemannian manifold $(M,g),$ one introduces the second order
differential
operators
$\Delta_q: \Omega^q(M) \to \Omega^{q}(M),$ the Laplace Beltrami
operator,
and $\Delta_q(t): \Omega^q(M) \to \Omega^{q}(M),$ the Witten Laplacian
(for the
1-form
$\alpha$) by
$$\Delta_q:= \delta^{q+1} \cdot d^q + d^{q-1}\cdot \delta^q,$$ and
 $$\Delta_q(t):= \delta^{q+1}(t)\cdot d^q(t) + d^{q-1}(t)\cdot \delta^
 q(t).$$
Note that $\Delta_q(0)= \Delta_q.$ In view of
 (2.3) -(2.10) and (2.12) one verifies

$$\leqalignno{\Delta_q(t)= \Delta_q + t(L_{-grad_g\alpha} +L_ {-grad_g
\alpha}^{\sharp}) +
t^2 ||\alpha||^2 Id
&&(2.16)\cr}$$
where $grad_g \alpha$ is the unique vector field defined by
$(grad_g\alpha)^\sharp = \alpha.$ One verifies that
$L_{-grad_g\alpha} +L_ {-grad_g\alpha}^{\sharp}$ is a zeroth  order
differential
operator.

The operators $\Delta_q(t)$ are  elliptic, essentially selfadjoint, 
and positive, hence their spectra, $\text{spect} \Delta_q(t),$ are contained in
$[0,\infty).$ Further
$$ \ker \Delta_q(t)= \{ \omega \in \Omega^q(M) | d^q(t)=0,
\delta^q(t)=0 \}$$
If $\alpha= dh$ and 
if $0$ is an eigenvalue of $\Delta_q(0),$ then $0$ is an eigenvalue of
$\Delta_q(t)$ for all $t$ and with the same
multiplicity: this because $d^q(t) = e^{-th}\cdot d^q \cdot e^{th\cdot}$ and  
$\delta^q(t)= e^{-th}\cdot \delta^q \cdot e^{th\cdot}.$

\vskip .2in

\subhead {c) Spectral gap theorems}
\endsubhead
\vskip .2in

If $\alpha$ is a closed 1-form we write $Cr(\alpha)$ for the
set of zeros of $\alpha.$ This notation is justified because
in a neighborhood of any connected component of $Cr(\alpha),\ $
$\alpha= dh$
for some smooth function $h,$ unique up to an additive constant.

The pair $(\alpha, g)$ is called a Morse pair, resp. Morse-Bott
pair, resp.
G-Morse pair, resp. normal G-Morse pair
if $(h,g)$ satisfies C1, resp. C'1, resp. G-C1, resp. normal G-C1.

In this subsection we will study the spectrum of $\Delta_q(t)$ for the
Witten deformation associated with
$(\alpha, g)$ being a Morse pair, resp. a Morse Bott pair or a  G-Morse pair.

In the case of a Morse-Bott or a  G-Morse pair we consider an
additional complex, namely 
the complex of the critical
 sets
$$(C^\ast, d^\ast): = \bigoplus_{\Sigma} (\Omega^{\ast
 -i(\Sigma)}(\Sigma, o(\nu_-)),d^{\ast-i(\Sigma)}).$$ 
Note that $(C^\ast, d^\ast) \ne (C^\ast, \partial^\ast)$
In the case of a Morse pair  this complex is trivial as it is 
concentrated in degree zero.
The metric $g$ induces a Riemannian metric on
$Cr(\alpha)= \cup \Sigma,$ hence the
complex $(C^\ast, d^\ast)$ gives rise to 
Laplacians 
$$\Delta'_q: =
\bigoplus_{\Sigma}\Delta_{q-i(\Sigma)}(\Sigma).$$
In the case of a Morse-Bott pair $\Delta_q(
\Sigma)$ denotes the Laplacian on  $q-$forms with
coefficients the orientation bundle in $o(\nu_-),$ of $\nu_-.$
In the case of a $G-$Morse pair,$\Delta_q(\Sigma)$ denotes 
the Laplacian acting on $\Omega^q(G/H_{\Sigma}, \det (\rho_-))$
where $G/H_{\Sigma}$
is equipped with the Riemannian metric induced by a scalar product on
$\frak g$ which
invariant with respect to the adjoint
action of $H_{\Sigma}.$
Here $H_{\Sigma}$ is the isotropy group of a point $x\in \Sigma,$
which is independent of $x$ up to conjugacy.
If the $G-$triangulation is normal then $\rho_-$ is trivial and
$\Omega^q(G/H_{\Sigma}, \det(\rho_-))= \Omega^q(G/H_{\Sigma}).$

Note that $\Delta(G/H_{\Sigma})=\hat{\oplus} \Delta^{\xi}_q(G/H_{\Sigma})$ and
$\Delta^{\xi}_q(G/H_{\Sigma})$ is
an endomorphism in the finite
dimensional vector space $\Omega_\xi^q(G/H_{\Sigma}, o(\rho_-)),$ which can be
explicitly calculated by elementary representation theory.

The following result
is essentially due to E.Witten.

\proclaim {Theorem 2.1} Suppose that $(\alpha, g)$ is a Morse pair.
Then there exist 
constants

\noindent $C_1, C_2, C_3$ and $T_0$ depending on $(\alpha, g)$  so that for
any
$t> T_0$

1): $spect \Delta_q(t) \cap (C_1 e^{-C_2t}, C_3t )=\emptyset$ 

\noindent and

2) the number of eigenvalues of $\Delta_q(t)$
in the interval $[0,C_1 e^{-C_2t}]$ counted with their multiplicity is
equal to the number of zeros of $\alpha$  of index $q.$
\endproclaim

The above theorem states the existence of a gap in the spectrum of
$\Delta_q(t),$ namely the open interval $(C_1 e^{-C_2t}, C_3 t) ,$
which widens to $(0,\infty)$ when $t\to \infty.$

Clearly $C_1, C_2, C_3$ and  $T_0$ determine a constant $T>T_0$, so
that
for $t\geq T,$$\ 1\in (C_1e^{-C_2t}, C_3t)$
and therefore
$$spect \Delta_q(t) \cap [0, C_1 e^{-C_2 t}] = spect \Delta_q(t)
\cap [0, 1] $$ and
$$spect \Delta_q(t) \cap [C_3t, \infty) = spect \Delta_q(t) \cap [1,
\infty).$$

For $t>T$ we denote by $\Omega^q(M)_{sm}(t)$ the finite dimensional
subspace of $\Omega^q(M)$
generated by the $q-$eigenforms
of $\Delta_q(t)$ corresponding to the eigenvalues of $\Delta_q(t)$
smaller than $1.$ Note that $\Omega^q(M)_{sm}(t)$ is of dimension 
$m_q$ where $m_q$ is the number of critical points 
of index $q$ of the closed
1-form $\alpha.$

The theory of elliptic operators implies that these eigenforms which are 
a priori elements
in the
$L_2-$completion of $\Omega^q(M),$ are actually smooth, i.e. 
in $\Omega^q(M).$ Note
that
$d(t)(\Omega^q(M)_{sm}(t))\subset \Omega^{q+1}(M)_{sm}(t),$
so that $(\Omega^*(M)_{sm}(t), d^*(t))$ is a finite dimensional
cochain subcomplex
of $(\Omega^*(M), d^*(t))$ and $(e^{th} ( \Omega^*(M)_{sm}(t)), d^*)$
is a
finite dimensional subcomplex of $(\Omega^*(M), d^*).$ Clearly the
$L_2-$orthogonal complement of $\Omega^\ast(M)_{sm}(t)$ in $\Omega^\ast(M)$
is also a closed Frechet subcomplex $(\Omega^\ast (M)_{la}(t), d^\ast(t))$
of $(\Omega^\ast (M), d^\ast(t))$  and we have the following decomposition
$$\leqalignno{(\Omega^\ast(M), d^\ast(t))= (\Omega^\ast (M)_{sm}(t), d
^\ast(t))\oplus
(\Omega^\ast (M)_{la}(t), d^\ast(t))&&(2.17)\cr}$$
with $(\Omega^\ast (M)_{la}(t), d^\ast(t))$ acyclic.

Let us consider now the case of a Morse-Bott pair. Denote by
${\lambda'}_{q,1} \leq {\lambda'}_{q,2}\leq \cdots
\leq {\lambda'}_{q,r}\leq \cdots$
be the spectrum of $\Delta'_q,$ and by $\lambda'_q$ the first nonzero
eigenvalue of $\Delta'_q.$
The following result is due to Helffer [H] (cf also [P]).

\proclaim {Theorem 2.1'} Suppose that $(\alpha,g)$ is a Morse-Bott
pair and $r\geq 1,$ an integer.
Then there exist positive constants
$C_1, C_2, T_0,$ so that for any $t>0$ and $t\geq T$ and $1\leq q\leq n$

$|\lambda_{q, r}(t)- \lambda'_{q,r}|< C_1 t^{-C_2}$ where
$\lambda_{q, r}(t)$ is the $r-$th eigenvalue of
$\Delta_q(t).$ 
\endproclaim

In particular one can find $T \geq T_0$ so that
for  $\lambda'= \inf \{\lambda'_q\}$ and $t>T$

1)
$Spect \Delta_q(t) \cap (C_1t^{-C_2}, -C_1t^{-C_2} +\lambda')=
 \emptyset$

2) $\lambda' /2 \in (C_1t^{-C_2}, -C_1t^{-C_2} +\lambda') .$

Therefore one can again produce a decomposition of the form
$$\leqalignno{(\Omega^\ast(M), d^\ast(t))= (\Omega^\ast (M)_{sm}(t), d
^\ast(t))\oplus
(\Omega^\ast (M)_{la}(t), d^\ast(t))&&(2.17')\cr}$$
where $(\Omega^*(M)_{sm}(t), d^*(t))$ is a finite dimensional cochain
subcomplex
of $(\Omega^*, d^*(t))$ with $(\Omega^q(M)_{sm}(t)$ given by the span of
the eigenforms
corresponding to eigenvalues
of $\Delta_q(t)$ smaller than $\lambda'$ and $(\Omega^*(M)(t)_{la}$
its orthogonal complement in
$\Omega^\ast (M).$ Clearly
$(\Omega^\ast (M)_{la}(t), d^\ast(t))$ is acyclic.

In [H] one finds $C_2 \geq 5/2$ and in [P]
$C_2> 1/2.$ 
Under additional hypothesis better estimate that the 
one stated in Theorem
2.1' can be obtained. As an example we mention 
the case of a $G-$Morse pair considered below.

Suppose $\mu:G\times M\to M$ is a smooth $G-$manifold, $G$ being a compact
Lie group,
$(\alpha, g)$ a $G-$Morse pair and $\xi$ an irreducible representation
of $G.$
Since $G$ is compact and the Riemannian metric $g$ is $G-$invariant 
$\Delta'$ and $\Delta_q(t)$ decompose orthogonally as 
$\Delta'_q= \hat{\sum}{\Delta'}_q^\xi$ and $\Delta_q(t)= \hat{\sum}\Delta_q^{\xi}(t)$

\noindent Let $\lambda^\xi_{q,1},\cdots \lambda^\xi_{q,N}(t)$ be the eigenvalues of
$({\Delta'}_q^\xi$ and $\lambda^\xi_{q,1}(t)\leq \cdots $
the ones of the Witten Laplacian $\Delta_q(t)^{\xi}.$ 

\vskip .2in

\proclaim {Theorem 2.1"} There exist the positive constants
$C_1, C_2, C_3$ and $T_0$ depending on $M, \alpha, g $ and $\xi$  so
that for any
$t> T_0$

1)
$\text{Spect} {\Delta_q(t)}^\xi \subset\bigcup_{i=1,\cdots,N}(-C_1
 e^{-tC_2}+\lambda^\xi_
{q,i},\
C_1 e^{-tC_2}+\lambda^\xi_{q,i})\cup
[C_3 t,\infty)$ and

2) the number of the eigenvalues of ${\Delta_q(t)}^\xi$
in the interval

\noindent $(-C_1 e^{-tC_2}+\lambda^\xi_{q,i},\
C_1 e^{-tC_2}+\lambda^\xi_{q,i})$  equals the multiplicity of
$\lambda^\xi_{q,i}$
\endproclaim

In particular one can canonically decompose
$(\Omega^*(M)_\xi , d^*(t)_\xi)$
in an orthogonal sum
$$\leqalignno{(\Omega^\ast_\xi(M), d^\ast(t)_\xi )= 
((\Omega^\ast_\xi(M))_{sm}(t), d^\ast(t)_\xi)\oplus
((\Omega^\ast_\xi(M))_{la}(t), d^\ast(t)_{\xi})&&(2.17")\cr}$$
where $((\Omega^q_\xi(M)(t))_{sm}(t))$ is the span of the eigenforms
corresponding to the eigenvalues
$\lambda^\xi_{q,i}$ and  $((\Omega^{\ast}_\xi(M))_{la}(t))$ is the  orthogonal
complement.
Note  that 

\noindent $((\Omega^{\ast}_\xi(M))_{la}(t)), d^\ast(t))$ is acyclic.

The proof of Theorem 2.1 will be given in subsection e) (cf [BZ1] and [BFKM]. The proof of
Theorem 2.1' and 2.1" can be
found in literature in [H] and [BFK5].

\vskip .2in

\subhead {d) Applications}
\endsubhead
\vskip .2in

{\it Morse inequalities:}

\vskip .1in

Let $\alpha$ be a closed 1-form and denote by $[\alpha]$ its
cohomology class,
which can be interpreted as a 1-dimensional representation of the
fundamental group of $M.$
Let $\beta_i (M,[\alpha]):= \dim H^i(M; [\alpha]).$
It is not hard to show that the integer valued function
$\beta(M,[t\alpha])$ is constant in $t$ for $t$ large
enough so that
$\hat{\beta}_i(M,[\alpha]):=\lim _{t\to \infty} \beta_i (M,[t\alpha])$
is well defined.

As an immediate consequence of the  decompositions discussed in section c)
we have
the following result: 
\proclaim
{Theorem 2.2}  Suppose that $\alpha$ is a Morse one form
and let $C_i:=\sharp (Cr_i
(\alpha)).$
Then for any integer $N$, we have
$$(-1)^N\sum_{i=0}^{N} (-1)^i C_i\geq (-1)^N\sum_{i=0}^{N}
(-1)^i {\beta}_i(M,[t \alpha])$$
for $t\geq T_0.$
If $\alpha$ is  a Morse Bott one form and $Cr(\alpha)$ is the union of
the closed connected submanifolds $\Sigma,$
then the above formula holds with $C_i$ given by
$C_i= \sum_{\Sigma} \dim H^{i-i(\Sigma)}(\Sigma, o(\nu_-)).$
\endproclaim
Clearly the inequalities remain true with ${\beta}_i(M,[t \alpha])$
replaced by $\hat{\beta}_i(M,[\alpha]).$ The above result is known as
the Morse inequalities when $\alpha= dh$ and
as the 
Novikov-Morse inequalities  when $\alpha$ is a closed 1-form.
It has a number of pleasant consequences in symplectic topology which
will be discussed below.
\vskip .2in

{\it Symplectic vector fields:}
\vskip .1in
Let $(M^{2n},\omega)$ be a symplectic manifold. The nondegenerated 2-
form $\omega$ establishes a bijective
correspondence $X\to \alpha_X$ between the set of smooth vector fields
$X$ on $M$ and the smooth 1-forms
$\alpha $ with $\alpha_X$ defined by the formula
$\alpha_X(Y)= \omega(X,Y)$ for any vector field $Y.$

Given a vector field $X$ consider
$\text{Zeros}(X):=\{x\in M |X(x)=0 \}.$
We say that $x^0\in \text{Zeros}(X)$ is nondegenerated if in one (and
then any)
coordinate system $\{x_1,\cdots x_{2n}\}$ around $x^0\in M$ the vector
field
$X=\sum_{i=1,\cdots,2n}a_i(x_1,\cdots, x_{2n})\partial/ \partial x_i$
satisfies
$\det(\partial a_i/\partial x_j(x^0_1,\cdots, x^0_{2n}))\ne 0$ where
$(x^0_1,\cdots, x^0_{2n})$ are the coordinates of $x^0.$
In this case we set
$$\leqalignno{ I(x^0):= \text{sign} \det(\partial a_i/\partial
x_j(x^0_1,\cdots, x^0_{2n})).&& (2.18)\cr} $$
The famous Hopf Theorem states:

\proclaim {Theorem} If $M$ is closed and all zeros of $X$ are
nondegenerated, then
$$\sum_{x\in \text{Zeros}(X)} I(x)= \chi (M)$$ where $\chi(M)$ denotes
the Euler Poincar\'e
characteristic of $M.$
In particular
$$\sharp (\text{Zeros}(X))\geq |\chi(M)|.$$
\endproclaim

\noindent This is the best general result about such vector fields. The estimate
is sharp.

The vector field $X$ is called symplectic if the 1-parameter (local)
group of
diffeomorphisms  induced by $X$ preserves the form $\omega,$ equivalently if
$L_X(\omega)=0$ or,
equivalently if $d(\alpha_X)=0.$
Clearly, $\text{Zeros}(X)= Cr(\alpha_X)$ and if all zeros of $X$ are
nondegenerated,
$\alpha_X$
is a Morse form and one can easily verify that $I(x)= (-1)^{i(x)}$
where $i(x)$
denotes the index of the critical point $x$ of the form $\alpha_X$.
Theorem 2.2 provides a considerable improvement of the Hopf Theorem, in
particular it says 
that for a symplectic vector
field with all zeroes nondegenerated
$$\leqalignno{\sharp (\text{Zeros}(X))\geq \sum_i \hat{\beta}_i(M,[
\alpha]).&& (2.19)\cr}$$
In fact, as shown in [BH2], it is possible to prove the existence of zeros
of a symplectic vector field in certain cases 
by "torsion methods"
even when $\hat{\beta}_i(M,[\alpha])=0$ for any $i.$

One can also apply the above theory to the study of some Lagrangian
intersections. A precise
situation
is the case of the symplectic manifold $T^\ast M,$ the cotangent
bundle of a smooth manifold $M,$
equipped with the canonical symplectic structure.

We are interested in the
intersection of
the zero section, the canonical Lagrangian in $T^\ast M,$ with the
image of a Lagrangian immersion
$i:\Sigma \to T^\ast M.$
In case that $\Sigma$ has a generating function $h:E\to R,$ where $E
$ is the total space of a smooth vector
bundle on $M,$ and
all critical points of $h$ are nondegenerated (cf [MS] for the definition of
a generating function), the count of the intersection points of
$i(\Sigma)$ and $M$
reduces to the count of zeroes of the closed form $i^\ast(\sigma)$ on
$\Sigma.$ As all zeros
are nondegenerated and Theorem 2.2 applies.
Here $\sigma$ is the canonical
1-form on $T^{\ast} M$ whose differential
$d(\sigma)$ defines the canonical symplectic structure on
$T^{\ast} M.$

\vskip .2in
\subhead {e) Sketch of  the proof of Theorems 2.1}
\endsubhead
\vskip .2in

The proof of Theorems 2.1 stated in the next section c) is based on a mini-max
criterion for detecting a
gap
in the spectrum of a positive selfadjoint operator in a Hilbert space
$H$ (cf 
 Lemma 2.3 below) and uses the  explicit formula for $\Delta_q(t)$ in
 admissible
 coordinates in a neighborhood of the set of critical points.

 \proclaim {Lemma 2.3} Let $A: H \to H$ be a densely defined (not
 necessary bounded )
 self adjoint positive operator in a Hilbert space $(H,<,>)$
and  $a,b$
two real numbers so that $0<a <b < \infty.$  Suppose that there exist
two closed subspaces $H_1$ and $H_2$ of $H$ with $H_1\cap H_2 =0$
and $H_1 + H_2 =H$ such that
 \newline (1) $<Ax_1, x_1> \ \leq \ a||x_1||^2$ for any $x_1\in H_1,$
 \newline (2) $<Ax_2, x_2> \ \geq \ b||x_2||^2$ for any $x_2\in H_2.$

 Then $spect A\bigcap (a,b)= \emptyset.$

 \endproclaim

 The proof of this Lemma is elementary (cf [BFK3] Lemma 1.2) and might be 
 a good exercise for
 the reader.
\vskip .1in

Consider $x\in Cr(\alpha)$ and choose admissible coordinates
$ (x_1,x_2, ...,x_n) $
in a neighborhood of $x$. With respect to these coordinates
$\alpha= dh,$
$$h(x_1,x_2, ...,x_n) =  - 1/2 (x_1^2+\cdots+ x_k^2) +1/2 (x_{k+1}^2+
\cdots +
x_n^2)$$
and  $g_{ij}(x_1,x_2, ...,x_n) = \delta_{ij},$ and hence by (2.16)
the operator $\Delta_q(t)$
 has the form
$$\leqalignno{\Delta_{q,k}(t)= \Delta_q + t M_{q,k} +t^2 (x_1^2 +
\cdots +x_n^2)Id
&&(2.20)\cr}$$
with $$\Delta_q (\sum _ I a_I(x_1,x_2, ...,x_n) dx_I)
= \sum_I (-\sum _{i=1}^n \frac{\partial^2}{\partial x_i^2} a_I (x_1,
x_2, ...,x_n))
dx_I,$$
and $M_{q,k}$ is the linear operator determined by
$$
\leqalignno{ M_{q,k} (\sum _ I a_I(x_1,x_2, ...,x_n) dx_I)=
\sum _I \epsilon_I^{q,k} a_I(x_1,x_2, ...,x_n) dx_I .&&(2.21)\cr}$$
Here  $I =(i_1,i_2\cdots i_q),$  $1\leq i_1 < i_2 \cdots< i_q \leq n,$
 $ dx_I= dx_{i_1}\wedge \cdots \wedge d_{i_q} $ and
$$\epsilon _I^{q,k} = -n+2k-2q +4 \sharp\{j|k+1 \leq i_j \leq n\},$$
where $\sharp A$ denotes the  cardinality of the set $A.$
Note that $\epsilon^{q,k}_I \geq -n$ and is $=-n$ iff $q=k.$

Let $\Cal S^q(\Bbb R^n)$ denote the space
of smooth $q-$forms
$\omega = \sum_I a_I(x_1,x_2, ...,x_n) dx_I$
with $a_I(x_1,x_2, ...,x_n)$ rapidly
decaying functions. The operator $\Delta_{q,k}(t)$ acting on
$\Cal S^q(\Bbb R^n)$ is globally elliptic (in the sense of [Sh1]
or [H\"o]), selfadjoint and positive. This operator is the harmonic
oscillator in $n$ variables acting on $q-$forms and its properties can
be derived
from the harmonic oscillator in one variable
$-\frac{d^2}{dx^2}+a +bx^2$
acting on functions. In particular the following result holds.

\proclaim{ Proposition 2.4} (1) $\Delta_{q,k} (t),$ regarded as an
unbounded
densely defined operator on the $L_2-$completion of
$\Cal S^q(\Bbb R^n),$ is
selfadjoint, positive and its spectrum is contained in
$2t\Bbb Z_{\geq 0}$
(i.e positive integer multiples of $2t$).

(2) $\ker \Delta_{q,k}(t)=0\  \text{if} \
k\ne q$ and $dim \ker\Delta_{q,q}(t) = 1.$

(3)
$\omega_{q,t}= (t/{\pi})^{n/4} e^{-t\sum_i x_i^2/2} dx_1\wedge\cdots
 \wedge dx_q$
is the generator of $\ker\Delta_{q,q}(t)$ with $L_2-$norm $1.$
\endproclaim

For a proof consult [BFKM] page 805.

Choose a smooth  function $\gamma_\eta(u),$
$\eta\in (0,\infty), \ u\in \Bbb R,$
which satisfies 
$$ \leqalignno{\gamma_\eta(u)=\quad
\left\{\aligned 1\ \text {if} \ & u\leq \eta/2 \\ 0 \ \text{if}\  &u
>\eta\endaligned \right\}.&&(2.23)\cr}$$

Introduce
$\tilde \omega_{q,t}^\eta \in \Omega^q_c(\Bbb R^n) $ defined by

$$\tilde\omega_{q,t}^\eta (x)= \beta_q(t)^{-1} \ \gamma_\eta (|x|)
\omega_{q,t}(x)$$ with $|x|= \sqrt{\sum_i x_i^2}$ and

$$\leqalignno{ \beta_q(t)= (t/{\pi})^{n/4}( \int_{\Bbb R^n}
\gamma_\eta^2(|x|) e^{-t \sum_ix_i^2} dx_1\cdots dx_n)^{1/2}.
&&(2.24)\cr}$$

The smooth form $\tilde\omega_{q,t}^\eta$ has its support in the ball
$\{|x|\leq \eta\},$ agrees with $\omega_{q,t}$ on the ball $\{|x|\leq \eta/2\}$ 
and satisfies

$$\leqalignno{< \tilde\omega_q^\eta(t), \tilde\omega_q^\eta(t)> =1
&&(2.25)\cr}$$ with respect to the scalar product $<.,.>$ on
$\Cal S^q(\Bbb R^n),$
induced by the Euclidean metric. The following proposition can be
obtained
by elementary calculations in coordinates in view of the explicit
formula of
$\Delta_{q,k}(t)$ (cf [BFKM], Appendix 2).

\proclaim{Proposition 2.5} For a fixed $r\in \Bbb N_{\geq 0}$ there
exist positive constants $C,C',C'', T_0,$ and $\epsilon_0$ so that $t >T_0$ and
$\epsilon <\epsilon_0$ imply

(1) $ |\frac{\partial^{|\alpha|}}{\partial x_1^{\alpha_1}\cdots
\partial x_n^{\alpha_n}}\Delta_{q,q}(t)
\tilde\omega_{q,t}^\epsilon(x)|\leq Ce^{-C't}$ for any
$x\in \Bbb R^n$ and multiindex $\alpha= (\alpha_1,\cdots,\alpha_n)$,
with
$|\alpha|=\alpha_1+\cdots+ \alpha_n \leq  r.$

(2) $<\Delta_{q,k}(t)\tilde\omega_{q,t}^{\epsilon},
\tilde \omega_{q,t}^{\epsilon}>\  \geq\  2t|q-k|$

(3) If $\omega \perp \tilde\omega^\epsilon _{q,t}$  with respect to
the scalar
product $<.,.>$
then $$<\Delta_{q,q}\omega,\omega>\  \geq \ C'' t||\omega||^2.$$

\endproclaim

For the proof of Theorem 2.1 (and  of Theorem 3.1 in the next lecture) we set the 
following notations.
We choose
$\epsilon >0$ so that for each $y\in Cr(h)$ there exists an
admissible coordinate chart
$\varphi_y: (U_y, y) \to (D_{2\epsilon},0) $ so that
$U_y\cap U_z=\emptyset$ for $y\ne z,$ $y,z\in Cr(h).$

Choose once and for all such an admissible coordinate chart for each
$y\in Cr_q(h).$
Introduce the smooth forms $\overline{\omega}_{y,t}\in \Omega^q(M)$
defined by
$$\leqalignno{ \overline{\omega}_{y,t} |_{M\setminus \varphi_y^{-1}
(D_{2\epsilon})}:=0, \ \
\overline{\omega}_{y,t}|_{\varphi_y^{-1}(D_{2\epsilon})}:=
\varphi_y^*(\tilde\omega^\epsilon_{q,t}).
&&(2.26)\cr}$$

For any given $t> 0$ the forms $\overline{\omega}_{y,t}\in\Omega^q(M), \ y\in Cr_q(h),$ are
orthonormal. Indeed,
if $y,z\in Cr_q(h),$ $y\ne z$ \  $\overline{\omega}_{y,t}$ and
$\overline{\omega}_{z,t}$
have disjoint support, hence are orthogonal, and because the support
of $\overline{\omega}_{y,t}$
is contained in an admissible chart,
$< \overline{\omega}_{y,t}, \overline{\omega}_{y,t}> = 1$ by (2.25).

For $t>T_0,$ with $T_0$ given by Proposition 2.5, we introduce
$J^q(t): C^q(X,\tau) \to \Omega^q(M)$
to be the linear map defined by
$$J^q(t)(E_y) = \overline{\omega}_{y,t}\ ,$$ where
$E_y\in C^q(X,\tau)$ is given by $E_y(z)= \delta_{yz}$ for
$y,z\in Cr(h)_q.$
$J_q(t)$ is an isometry, thus in particular injective.

{\bf Proof of Theorems 2.1:} (sketch). Take $H$ to be the
$L_2-$completion of
$\Omega^q(M)$
with respect to the scalar product $<., .>,$
$H_1:= J^q(t)(C^q(M,\tau))$ and
$H_2= H_1^{\perp}.$
Let $T_0, C, C', C''$ be given by Proposition 2.5 and define
$$C_1:= \inf _{z\in M'} || grad_g \alpha(z)||,$$ with
$\ M'= M\setminus \bigcup_{y\in Cr_q(\alpha)}
 \varphi_y^{-1}(D_\epsilon),$and
$$C_2= \sup_{x\in M} ||(L_{ -grad_g \alpha} +L_ {-grad_g \alpha}^{
\sharp})(z)||.$$
Here $||grad_g \alpha(z) ||$ resp.
$||(L_{-grad_g \alpha} +L_ {-grad_g \alpha}^{\sharp})(z)||$ denotes
the norm of
the vector $ grad_g \alpha(z)\in T_z(M)$ resp. of the linear map
$(L_{-grad_g \alpha} +L_ {-grad_g \alpha}^{\sharp})(z):
 \Lambda^q(T_z(M)) \to
 \Lambda^q(T_z(M))$
with respect to the scalar product  induced in $T_z(M)$ and
$\Lambda^q(T_z(M))$ by $g(z).$
Recall that if $X$ is a vector field then $L_X +L_X^{\sharp}$ is a
zeroth
order differential operator,
hence an endomorphism of the bundle $\Lambda^q(T^*M)\to M.$

We can use the constants
$T_0, C, C', C'', C_1, C_2$ to construct $C'''$ and $\epsilon_1$ so
that for $t> T_0$
and $\epsilon <\epsilon_1,$
we have
$< \Delta_q(t) \omega, \omega > \geq C_3t <\omega, \omega>$ for any
$\omega \in H_2$
 (cf. [BFKM], page 808-810).

Now one can apply Lemma 2.3 whose hypotheses are satisfied for
$a= Ce^{-C't}, b =C''' t$ and $t>T_0.$
This concludes the first  part of Theorem 2.1.

 Let $Q_q(t),$  $t>T_0$ denote the orthogonal projection in $H$ on the
 span of  the
 eigenvectors
 corresponding the eigenvalues smaller than $1.$ In view of the
 ellipticity of
 $\Delta_q(t)$ all
 these eigenvectors are smooth $q-$forms. An additional important
 estimate is
given by the following Proposition:

\proclaim{Proposition 2.6} For $r\in \Bbb N_{\geq 0}$ one can find
$\epsilon_0 >0 $ and $C_3, C_4$ so that
for $t> T_0$ as constructed above, and any $\epsilon <\epsilon_0$
one has, for any $v\in C^q(M,\tau)$ 
$$(Q_q(t) J^q(t)- J^q(t))(v) \in \Omega^q(M),$$
and for $0\leq p\leq r,$
$$ ||(Q_q(t) J^q(t) - J^q(t)) (v)||_{C^p} \leq C_3 e^{- C_4t}
||v||,$$
where $||\cdot||_{C^p}$ denotes the 
$C^p-$norm. 
\endproclaim

The proof of this Proposition is contained in [BZ1], page 128 and [
BFKM] page 811.
Its proof requires (2.16), Proposition 2.5 and general estimates
coming from the
ellipticity of $\Delta_q(t).$

Proposition 2.6 implies that for $t$ large enough, say $t> t_0,$
$\Cal I^q(t) := Q_q(t) J^q(t)$ is bijective, which  finishes the proof
of
Theorem 2.1.

\vskip .5in

{\bf Lecture 3: Helffer Sj\"ostrand Theorem, an asymptotic improvement of the Hodge-de
Rham
 theorem }

In this section we will formulate the result of Helffer Sj\"ostrand
for a
generalized triangulation and
its analogue for a $G-$generalized triangulation.

\vskip .2in
 \subhead {a. Hodge-de Rham  theorem and its (asymptotic)
 generalization}
 \endsubhead

First let us recall the classical Hodge-de Rham  theorem.

\noindent Consider a closed Riemannian manifold $(M,g)$ and a
simplicial smooth
triangulation $\tau.$
The Riemannian metric $g$  provides a scalar product in the Frechet
space of smooth forms
and then the Laplace Beltrami operators
$\Delta_i:\Omega^i(M) \to \Omega^i(M).$
Let us denote by $\Cal H^i:= \ker \Delta_i.$
The simplicial triangulation $\tau$ provides the (finite dimensional)
cohomology
vector spaces $H^\ast_{\tau}(M).$

\proclaim
{Theorem } (Hodge-de Rham )

Given a Riemannian manifold $(M,g)$ equipped with a smooth
triangulation $\tau$
one can produce:

1) a canonical orthogonal decomposition
$\Omega^\ast(M)= \Cal H^\ast\oplus
\Omega^{\ast}_2(M),$ with $\Cal H^\ast=
Ker \Delta_\ast$ a finite dimensional graded vector space
and
$\Omega^{\ast}_2(M)= d(\Omega^{\ast -1}(M))\oplus
 d^{\sharp}(\Omega^{\ast
+1}(M))$

2) a canonical linear isomorphism (whose inverse is induced from
integration
 of forms on simplexes)
$J: H^{\ast}_{\tau}(M,\Bbb R)\to \Cal H^\ast.$
\endproclaim

The above theorem provides a canonical realization of the cohomology,
calculated with the help of the smooth
 triangulation $\tau,$
as differential forms, harmonic with respect to the given Riemannian
metric $g$.

One can improve the above result by realizing the full geometric
complex defined by
the triangulation as a subcomplex  of differential forms but the "
canonicity" statement
remains true only asymptotically.
More precisely one can show:
\vskip .2in

\proclaim {Theorem }
Given a Riemannian manifold $(M,g)$ and a smooth triangulation $\tau$
for
$t\in \Bbb R$ large enough
one can produce

1) a smooth one-parameter family of orthogonal decompositions
$$(\Omega^\ast (M), d)=(\Omega^\ast (M)_0(t), d)\oplus (\Omega^\ast
(M)_1(t), d)$$
with  $(\Omega^\ast (M)_0(t), d)$ a finite dimensional complex, which
is $O(1/t)$ canonical
 \footnote {i.e. for any two such possible decompositions the finite
 dimensional subspaces
 $\Omega^\ast (M)_0(t)$ are at an $O(1/t)$ distance with respect to
 the scalar product
 induced by the metric $g.$}, and 

2) a smooth family of isomorphisms
$\Cal I^\ast(t) :C^\ast (M,\tau) \to \Omega(M)_0^\ast
(t)$
so that the composition $\Cal I^\ast(t)\cdot S^\ast(t),$ where 
$S^\ast(t)$ is the scaling isomorphism
$$\leqalignno{S^\ast(t):(C^\ast(M,\tau),\delta^\ast_{\tau,o} ) \to (C^\ast (
M,\tau), \delta^\ast_{\tau,o}(t)) &&(3.1)\cr}$$
defined by $S^q(t)(E_x) = (\frac{\pi}{t})^{(n-2q)/4} e^{-th(x)} E_x, x\in Cr_q(h),$ is of the
form $I^\ast+O(1/t)$ with $I^\ast $ an isometry.
Moreover $\Cal I^\ast(t)$ is $O(1/t)$ canonical
\footnote {i.e. for any two such possible $\Cal I(t)'$s,
$\Cal I_1(t)\  \text {and} \ \Cal I_2(t),$
$||\Cal I_1(t) -\Cal I_2(t)|| = O(1/t)$}.
\endproclaim

In view of a result of Pozniak which claims that any smooth
triangulation can be
realized as a generalized  triangulation (cf section 1 a), O.1) the
above theorem
is a straightforward reformulation of
Theorem 3.1 below, proven by Helffer and Sj\"ostrand.

\vskip .1in
\subhead {b. Helffer Sj\"ostrand theorem (Theorem 3.1)}
\endsubhead

We consider only the case of a generalized triangulation  and of a G-
generalized triangulation.
We pick up orientations $o$ as indicated in section 1
and consider the scaling (3.1)
$S^q(t): (C^\ast (M,\tau),\delta^\ast_{\tau,o}) \to (C^\ast
 (M,\tau),\delta^\ast_{\tau, o}(t))$
and, for $t$ large enough, the compositions $L(t)$ and
$L(t)_\xi$ defined by the following diagram

$$
\CD
(\Omega^{\ast}(M),d^\ast(t))@>e^{th}>>(\Omega^{\ast}(M),d^\ast)@>Int^
\ast>>(C^{\ast}(M,\tau),\delta ^\ast_{\tau, o})
\\
 @AinAA                                           @.
 @VS^{\ast}(t)VV                         \\
(\Omega^{\ast}_{sm}(M)(t),d^\ast(t))  @.             @.             (C^{
\ast}(M,\tau),\delta^\ast_{\tau, o}(t))
\endCD
$$
and
$$
\CD
(\Omega^{\ast}(M)_\xi,(d^\ast(t))_\xi)@>e^{th}>>(\Omega^{\ast}(M)_\xi,
d^\ast_\xi)@>Int^\ast>>(C^{\ast}(M,\tau)_\xi,
\delta^\ast_\xi)
\\
 @AinAA                                           @.
 @VS^{\ast}(t)VV                         \\
((\Omega^{\ast}_{sm}(M)(t))_\xi,d^\ast(t)_\xi)  @.             @.
(C^{\ast}(M,\tau)_\xi,\delta^\ast_\xi(t))
\endCD
$$
The following theorems are reformulations of a theorem due to Helffer-
Sj\"ostrand,
cf [HS2].

\proclaim
{Theorem 3.1} (Helffer-Sj\"ostrand). Given $M$ a
closed manifold and $\tau=(g,h)$ a generalized triangulation, there
exists $T >0,$ depending
on $\tau,$ so that for $t> T,$
$L^\ast (t)$ is an isomorphism of cochain
complexes.
Moreover, there exists a family of isometries
$R^q(t): C^q(M,\tau) \to \Omega^q(M)_{sm}(t)$ of finite dimensional
vector
spaces so that
$ L^q(t)\cdot R^q(t)  = Id + O(1/t).$
\endproclaim

and
\proclaim
{Theorem 3.1"} Given a
closed $G-$manifold $M,$ with $G$ being a compact Lie group
and $\tau=(g,h)$ a $G-$generalized triangulation and an 
irreducible
representation $\xi,$ there
exists $T >0,$ depending
on $\tau$ and $\xi,$ so that for $t> T,$ $L^\ast(t)_\xi$ is an
isomorphism of cochain
complexes.
Moreover, there exists a family of isometrics
$R^q(t)_{\xi}: C^q(M,\tau)_\xi \to (\Omega^q(M)_{sm}(t))_\xi$ of
finite dimensional vector
spaces so that
$ L^q(t)_{\xi}\cdot R^q(t)_{\xi}  = Id + O(1/t).$
\endproclaim

It is understood that
$C^q(M,\tau)$ is equipped with
the canonical scalar product defined in section 1,
and $\Omega^q(M)_{sm}(t)$ with the scalar product $<\ ,\ >$ defined
by (2.15).

One can also prove an analogue of Theorem 3.1 for MB-generalized
triangulations. This requires a finite dimensional cell complex instead
of $(C^\ast, \partial^\ast)$.
Such complex is described in [BH] where an analogue of Theorem 3.1 for
a Morse-Bott form is
established.

{\bf Sketch for the proof of Theorems 3.1:}
The proof is a continuation of the proof of Theorem 2.1 in the
previous lecture and we use the same notation.
Therefore we invite the reader to review the section e) of Lecture 2 .

Let $T_0$ be provided by Proposition 2.6. For $t\geq T_0,$ let
$R^q(t)$ be
the isometry defined by
$$\leqalignno{ R^q(t):= J^q(t) (J^q(t)^{\sharp} J^q(t))
^{-1/2}
&&(3.2)\cr}$$
and introduce 
$U_{y,t}:= R^q(t) (E_y) \in \Omega^q(M)$ for any $y\in Cr_q(h)= Cr_q(dh).$
Proposition 2.6 implies that there
exists $\epsilon >0,$  $t_0$ and
$C$ so that for any $t > t_0$  and any $y\in Cr(h)_q$ one has

$$\leqalignno{ \sup_{z\in M\setminus \varphi_y^{-1}(D_\epsilon)}
||U_{y,t}(z)|| \leq Ce^{-\epsilon t} &&(3.3)\cr}$$
and
$$\leqalignno{ ||U_{y,t}(z) - \overline\omega_{y,t}(z)||
\leq C \frac{1}{t},\  \text{for any } \  z\in W_y^-\cap \varphi_y^{-1}(D_
\epsilon).
&&(3.4)\cr}$$

To check Theorem 3.1 it suffices to show that

$$| \int_{W_{x'}^-} U_{x,t} e^{th}\  -\ (\frac{t}{\pi})^{\frac{n-2q}{4}} e
^{-th(x)}
\delta_{xx'} |
\leq C''\ \frac{1}{t} $$
for some $C''> 0$ and any $x,x'\in  Cr(h)_q.$

If $x\ne x'$ this follows from (3.3). If $x=x'$ from (3.3) and (3.4).
\hfill $\square$

\proclaim{5. Extensions and a survey of other applications}  \endproclaim

1. One can relax the definition of admissible coordinates in C1, and
C1" by
dropping the requirements on the metric. Theorems 2.1, 2.1" and 3.1 and
3.1"
remain true as stated; however almost all calculations
will be longer since the explicit formulae for $\Delta_q(t)$ and its
spectrum
when regarded on $\Cal S^*(\Bbb R^n)$ will be more complicated.

2. One can provide an analogue of Theorems 2.1  and 3.1 in the case of
a closed
one form. This will be elaborated in a forthcoming paper [BH2].

3. One can twist both complexes $(\Omega^*, d^*)$ and
$(C^*(M,\tau), \partial^*)$ resp
$(C^*(M,\tau),D^\ast)$
by a finite dimensional representation of the fundamental group,
$\rho: \pi_1(M) \to GL(V).$
In this case additional data is necessary: a Hermitian structure $\mu$
on the flat bundle $\xi_\rho$ induced by
$\rho.$
The "canonical" scalar product on $(C^*(M,\tau,\rho), \partial^*),$
in the case of a generalized triangulation will be obtained by using
the critical points (the cells of the
generalized triangulation)
and the Hermitian scalar product provided by $\mu$ in the fibers
of $\xi_\rho$ above the critical points. The  de-Rham complex in this
case is
replaced by
$(\Omega^*(M,\rho), d_\rho^*)$ the  de-Rham complex of differential
forms with coefficients in
$\xi_\rho$ whose  differential is given by the covariant
differentiation w.r. to the canonical flat
connection in $\xi_\rho.$ To have a scalar product on the spaces of
smooth forms,
in addition to the Riemannian metric on $M$ one needs
a Hermitian structure $\mu$ (cf. [BFK1] or [BFK4]) in $\xi_\rho.$ The
statements of Theorems 2.1
and 3.1 remain the same.
Under the hypotheses that the Hermitian structure is parallel in small
neighborhoods of the critical points, the  proofs
remain the same.
An easy continuity argument permits to reduce the case of an arbitrary
Hermitian structure
to the previous one, by taking $C^0$ approximation of a given
Hermitian
structure by Hermitian structures which are parallel near the critical
points.
Since the Witten Laplacians do not involve derivatives of the
Hermitian structure
such a reduction is possible.
If the representation is a unitary representation on a finite
dimensional Euclidean
space
one has a canonical Hermitian structure in $\xi_\rho$  which is parallel with
respect to the
flat
canonical connection in $\xi_\rho.$
This extension was used in the new proofs of the Cheeger- Muller
theorem
and its extension
about the  comparison of the analytic and the Reidemeister torsion,
cf. [BZ],
[BFK1], [BFKM], [BFK4].

4. One can further extend the WHS-theory to the case where  $\rho$ is a
special
type of an infinite
dimensional representation, a representation of the fundamental group
in an $\Cal A-$
Hilbert module of finite type. This extension was done in [BFKM] for
$\rho$
unitary and in [BFK4] for
$\rho$ arbitrary. In this case the Laplacian $\Delta_q(t)$ do
not have discrete spectrum
and it seems quite remarkable that Theorems 2.1 and 3.1 remain true.
It is even more
surprising that
exactly the same arguments as presented above can be adapted to prove
them.
A particularly interesting situation is the case of the left regular
representation
of a countable group
$\Gamma$ on the Hilbert space $L_2(\Gamma)$ when regarded as an
$\Cal N(\Gamma)$ right
Hilbert module of the von Neumann algebra $\Cal N(\Gamma),$
cf.[BFKM] for definitions.
One can prove that Farber extended $L_2-$cohomology of $M,$ a compact
smooth
manifold with infinite fundamental group defined analytically (i.e.
using
differential forms and a Riemannian metric) and combinatorially (i.e
using a
triangulation) are isomorphic and therefore
the classical
$L_2-$Betti numbers and Novikov-Shubin invariants defined
analytically and
combinatorially are the same. For this last fact 
see [BFKM] section 5.3.

The WHS-theory was a fundamental tool in the proof of the equality of
the
\newline $L_2-$analytic and the $L_2-$Reidemeister
torsion presented [BFKM].

5. One can further extend Theorems 1.1, 1,1', 2.2, 2.2', 3.3, 3.3', to bordisms
$(M,\partial_-M,
\partial_+M),$
and $\rho$ a representation of $\Gamma=\pi_1(M)$ on an $\Cal A-$Hilbert module of
finite type.
In this case one has first to extend the concept of generalized
triangulation  to
such bordisms. This will
involve a pair $(h, g)$ which in addition to the requirements C1-C3 is
supposed to satisfy the following assumptions:
$g$ is product like near $\partial M=\partial_-M \cup \partial_+M,$
the function $h:M \to [a, b]$ satisfies $h^{-1}(a)= \partial_-M,$
$h^{-1}(b)= \partial_+M,
\  a,b$ regular values, and is linear on the geodesics normal to
$\partial M$ near $\partial M.$ In case $h$ is replaced by a closed 1-form the
requirement
is that this form vanishes on $\partial M.$
This extension was partly done in [BFK2] and was used to prove gluing
formulae for
analytic torsion and to extend the results of [BFKM] to manifolds
with boundary.

5. One can actually extend the WHS-theory to the case where $h$ is
a generalized Morse function, i.e. the critical points are either
nondegenerated or birth-death. This extension is much more subtle and
very important. Beginning work in this direction was done by Hon Kit
Wai
in his OSU dissertation.

\vfill \eject

\end

\enddocument